\documentclass[twocolumn]{autart}    

\usepackage[square,comma,sort&compress]{natbib}
\usepackage[dvips,colorlinks,bookmarksopen,bookmarksnumbered,citecolor=blue,urlcolor=blue]{hyperref}

\usepackage{graphics}
\usepackage{amssymb}
\usepackage{amsmath}
\usepackage{color}
\usepackage{xspace}
\usepackage{caption}
\usepackage{algpseudocode}
\usepackage{bbm}
\usepackage{comment}
\usepackage{algorithmicx}
\usepackage{subfig}
\usepackage{psfrag}
\usepackage{tikz}
\usetikzlibrary{shapes,arrows}
\usetikzlibrary{positioning, decorations.markings}
\usetikzlibrary{calc}

\tikzstyle{block}=[draw opacity=0.7,line width=1.4cm]

\definecolor{CranJ}{cmyk}{0,0.69,0.54,0.04} 
\definecolor{PinkJ}{cmyk}{0,0.71,0.43,0.12} 
\definecolor{Cran}{cmyk}{0,0.73,0.41,0.29} 
\definecolor{VRed}{cmyk}{0,0.75,0.25,0.2} 
\definecolor{ORed}{cmyk}{0,0.75,0.75,0} 
\definecolor{CBlue}{cmyk}{1,0.25,0,0} 

\newcommand{\oprocendsymbol}{\hbox{$\bullet$}}
\newcommand{\oprocend}{\relax\ifmmode\else\unskip\hfill\fi\oprocendsymbol}
\def\eqoprocend{\tag*{$\bullet$}}

\newcommand{\longthmtitle}[1]{\mbox{}\textit{{(#1):}}}

\newcommand{\VV}{\mathcal{V}}
\newcommand{\EE}{\mathcal{E}}
\newcommand{\GG}{\mathcal{G}}
\newcommand{\PP}{\mathcal{P}}
\newcommand{\LL}{\vect{L}}
\newcommand{\lL}{\vect{\mathsf{L}}}
\newcommand{\rR}{\vect{\mathsf{R}}}
\newcommand{\PPi}{\vect{\Pi}}
\newcommand{\pPi}{\vect{\mathsf{\Pi}}}
\newcommand{\MT}{\overline{M}}
\newcommand{\mT}{\underline{m}}

\newcommand{\SLya}{\mathcal{S}_0}
\newcommand{\SBar}{\bar{\mathcal{S}}_0}

\newcommand{\real}{{\mathbb{R}}}
\newcommand{\reals}{{\mathbb{R}}}
\newcommand{\realpositive}{{\mathbb{R}}_{>0}}
\newcommand{\realnonnegative}{{\mathbb{R}}_{\ge 0}}

\newcommand{\integernonnegative}{{\mathbb{Z}}_{\ge0}}
\newcommand{\Lnorm}{\left\|}
\newcommand{\Rnorm}{\right\|}
\newcommand{\naturals}{{\mathbb{N}}}
\newcommand{\eps}{\epsilon}

\newcommand{\Hlambda}{\hat{\lambda}}

\newcommand{\argmax}{\operatorname{argmax}}

\newcommand{\map}[3]{#1:#2 \rightarrow #3}

\newcommand{\until}[1]{\in\VV}

\newcommand{\setdef}[2]{\{#1 \; |\; #2\}}
\newcommand{\vect}[1]{\boldsymbol{\mathbf{#1}}}
\newcommand{\Bvect}[1]{\bar{\boldsymbol{\mathbf{#1}}}}
\newcommand{\Tvect}[1]{\tilde{\boldsymbol{\mathbf{#1}}}}
\newcommand{\Hvect}[1]{\hat{\boldsymbol{\mathbf{#1}}}}

\newcommand{\dvect}[1]{\dot{\vect{#1}}}

\newcommand{\Sym}[1]{\operatorname{Sym}(#1)}
\newcommand{\Diag}[1]{\operatorname{Diag}(#1)}
\newcommand{\kronecker}{\raisebox{1pt}{\ensuremath{\:\otimes\:}}}

\newcommand{\sumSig}[1]{\sum\nolimits_{i=1}^N\!\!\!#1}

 \newcommand{\boxend}{\hfill \ensuremath{\Box}}

\makeatletter
\renewcommand*{\@opargbegintheorem}[3]{\trivlist
      \item[\hskip \labelsep{\emph{ #1\ #2}}] \emph{(#3):}\ \itshape}
\makeatother

\begin{document}

\begin{frontmatter}
   \runtitle{Distributed convex optimization}
  
   \title{Distributed convex optimization via continuous-time
     coordination algorithms with discrete-time
     communication\thanksref{footnoteinfo}} 
  
  \thanks[footnoteinfo]{Corresponding author: S.S. Kia}
  
  \author[Paestum]{Solmaz S. Kia}\ead{solmaz@uci.edu} \quad 
    \author[Rome]{Jorge Cort\'es}\ead{cortes@ucsd.edu}
  \quad 
  \author[Rome]{Sonia Mart{\'\i}nez}\ead{soniamd@ucsd.edu}
  
  \address[Paestum]{Department of Mechanical and Aerospace
    Engineering, University of California, Irvine}
  \address[Rome]{Department of Mechanical and Aerospace Engineering,
    University of California, San Diego}
  \begin{abstract}
    This paper proposes a novel class of distributed continuous-time
    coordination algorithms to solve network optimization problems
    whose cost function is a sum of local cost functions associated to
    the individual agents.  We establish the exponential convergence
    of the proposed algorithm under (i) strongly connected and
    weight-balanced digraph topologies when the local costs are
    strongly convex with globally Lipschitz gradients, and (ii)
    connected graph topologies when the local costs are strongly
    convex with locally Lipschitz gradients.  When the local cost
    functions are convex and the global cost function is strictly
    convex, we establish asymptotic convergence under connected graph
    topologies.  We also characterize the algorithm's correctness
    under time-varying interaction topologies and study its privacy
    preservation properties. Motivated by practical considerations, we
    analyze the algorithm implementation with discrete-time
    communication.  We provide an upper bound on the stepsize that
    guarantees exponential convergence over connected graphs for
    implementations with periodic communication.  Building on this
    result, we design a provably-correct centralized event-triggered
    communication scheme that is free of Zeno behavior. Finally, we
    develop a distributed, asynchronous event-triggered communication
    scheme that is also free of Zeno with asymptotic convergence
    guarantees.  Several simulations illustrate our results.
  \end{abstract}

\end{frontmatter}

\section{Introduction}

An important class of distributed convex optimization problems
consists of the (un-)constrained network optimization of a sum of
convex functions, each one representing a local cost known to an
individual agent. Examples include distributed parameter
estimation~\citep{SSR-AN-VVV:08b,PW-MDL:09}, statistical learning~\citep{SB-NP-EC-BP-JE:10}, and optimal resource
allocation over networks~\citep{RM-SL:06}.  To find
the network optimizers, we propose a coordination model where each
agent runs a purely local continuous-time evolution dynamics and
communicates at discrete instants with its neighbors.  We are
motivated by the desire of combining the conceptual ease of the
analysis of continuous-time dynamics and the practical constraints
imposed by real-time implementations. Our design is based on a novel
continuous-time distributed algorithm whose stability can be analyzed
through standard Lyapunov functions.

\emph{Literature review}: 
In distributed convex optimization, most coordination algorithms are
time-varying, consensus-based dynamics~\citep{AN-AO:09,BJ-MR-MJ:09,SB-NP-EC-BP-JE:10,MZ-SM:09c,JD-AA-MW:12}
implemented in discrete time. Recent work~\citep{JW-NE:11,BG-JC:14-tac,FZ-DV-AC-GP-LS:11,JL-CYT:12} has introduced
continuous-time dynamical solvers whose convergence properties can be
analyzed via classical stability analysis. This has the added
advantage of facilitating the characterization of properties such as
speed of convergence, disturbance rejection, and robustness to
uncertainty.~\cite{JW-NE:11} establish asymptotic convergence under
connected graphs and~\cite{BG-JC:14-tac} extend the design and
analysis to strongly connected, weight-balanced digraphs.
The continuous-time algorithms in~\citep{FZ-DV-AC-GP-LS:11,JL-CYT:12}
require twice-differentiable, strictly convex local cost functions to
make use of the inverse of their Hessian and need a careful
initialization to guarantee asymptotic convergence under undirected
connected graphs.  The novel class of continuous-time algorithms
proposed here (upon which our implementations with discrete-time
communication are built) do not suffer from the limitations discussed
above and have, under some regularity assumptions, exponential
convergence guarantees.  Our work here also touches, albeit slightly,
on the concept of privacy preservation, see
e.g.,~\citep{PCW-GA-CF-JSB:13,FY-SS-SV-YQ:13} for recent works in the
context of distributed multi-agent optimization. Regarding
event-triggered control of networked systems, recent years have seen
an increasing body of work that seeks to trade computation and
decision making for less communication, sensing or actuator effort,
see e.g.~\citep{WPMHH-KHJ-PT:12,MMJ-PT:11,XW-MDL:11}. Closest to the
problem considered here are works that study event-triggered
communication laws for average consensus, see
e.g.,~\citep{DVD-EF-KHJ:12,EG-YC-HY-PA-DC:13,CN-JC:14-auto}.  The
strategies proposed in~\citep{PW-MDL:09} save communication effort in
discrete-time implementations by using local triggering events but are
not guaranteed to avoid Zeno behavior, i.e., an infinite number of
triggered events in a finite period of time.  Our goal is to combine
the best of both approaches by synthesizing provably-correct
continuous-time distributed dynamical systems which only require
communication with neighbors at discrete instants of time.  We are
particularly interested in the opportunistic determination of this
communication times via event triggering schemes.

\emph{Statement of contributions}: We propose a novel class of
continuous-time, gradient-based distributed algorithms for network
optimization where the global objective function is the sum of local
cost functions, one per agent.  We prove that these algorithms
converge exponentially under strongly connected and weight-balanced
agent interactions when the local cost functions are strongly convex
and their gradients are globally Lipschitz. Under connected,
undirected graphs, we establish exponential convergence when the local
gradients are just locally Lipschitz, and asymptotic convergence when
the local cost functions are simply convex and the global cost
function is strictly convex.  We also study convergence under networks
with time-varying topologies and characterize the topological
requirements on the communication graph, algorithm parameters, and
initial conditions necessary for an agent to reconstruct the local
gradients of local cost functions of other agents.  
Our technical approach builds on the identification of strict Lyapunov functions.
The availability of these functions enable our ensuing design
of provably-correct continuous-time implementations with discrete-time
communication.  In particular, for networks with connected graph
topologies, we obtain an upper bound on the suitable stepsizes that
guarantee exponential convergence under periodic
communication. Building on this result, we design a centralized,
synchronous event-triggered communication scheme with an
exponential convergence guarantees and Zeno-free behavior. Finally,
we develop a Zeno-free asynchronous event-triggered communication
scheme whose execution only requires agents to interchange information
with their neighbors and establish its 
exponential convergence to a neighborhood of the network optimizer.
 Several simulations illustrate our results.

\section{Preliminaries}\label{sec::prelim}

In this section, we introduce our notation
and some basic concepts from convex functions and graph theory.
Let $\reals$ and $\integernonnegative$ denote, respectively, the set
of real and nonnegative integer numbers. We use $\Re(\cdot)$ to
represent the real part of a complex number. The transpose of a matrix
$\vect{A}$ is~$\vect{A}^\top$.  We let $\vect{1}_n$
(resp. $\vect{0}_{n}$) denote the vector of $n$ ones (resp. $n$
zeros), and denote by $\vect{\sf{I}}_n$ the $n\times n$ identity
matrix. We let $\pPi_n = \vect{\sf{I}}_n - \frac{1}{n}
\vect{1}_n\vect{1}_n^\top$.  When clear from the context, we do not
specify the matrix dimensions. For $ \vect{A}\!\in\real^{n\times m}$
and $\vect{B}\in\real^{p\times q}$, we let $
\vect{A}\!\kronecker\!\vect{B}$ denote their Kronecker product.  For
$\vect{u}\in\reals^d$, $\|\vect{u}\|=\sqrt{\vect{u}^\top\vect{u}}$
denotes the standard Euclidean norm.  For vectors
$\vect{u}_1,\cdots,\vect{u}_m$, $\vect{u} =
(\vect{u}_1,\cdots,\vect{u}_m)$ is the aggregated vector.  In~a
networked system, we distinguish the local variables of each agent by
a superscript, e.g., $\vect{x}^i$ is the state~of agent $i$. If
$\vect{p}^i\in\reals^d$ is a variable of agent $i$, the aggregated
$\vect{p}^i$'s of the network of $N$ agents is~$\vect{p} =
(\vect{p}^1,\cdots,\vect{p}^N) \in (\reals^d)^N$. For convenience, we
use $\vect{\mathsf{r}} \in \real^N$ and $\rR \in \real^{N \times (N-1)}$,
\vspace{-0.1in}
\begin{align}
  & \vect{\mathsf{r}}\!=\!\frac{1}{\sqrt{N}}\vect{1}_N,
  ~\vect{\mathsf{r}}^\top\rR \!=\!\vect{0},~ \rR^\top\rR \! =\!
  \vect{\mathsf{I}}_{N-1},~ \rR\rR^\top\!
  =\!\pPi_N. \label{eq::Rrs} \vspace{-0.1in}
\end{align}
A differentiable function $f: \reals^d\to\reals$ is \emph{strictly
  convex} over a convex set $C \subseteq \reals^d$~iff $
  (\vect{z}-\vect{x})^\top(
  \nabla f(\vect{z})-\nabla f(\vect{x}))>
  0$ for all $\vect{x},\vect{z}\in C$ and $\vect{x}\neq\vect{z}$,
and it is \emph{$m$-strongly convex} ($m>0$) iff
$(\vect{z}-\vect{x})^\top( \nabla f(\vect{z})-\nabla f(\vect{x}))\geq
m\|\vect{z}-\vect{x}\|^2$, for all $\vect{x},\vect{z}\in C$.
A function $\vect{f}: \reals^d\to\reals^d$ is Lipschitz with constant
$M>0$, or simply \emph{$M$-Lipschitz}, over a set $C\subset \reals^d$
iff $ \Lnorm \vect{f}(\vect{x})-\vect{f}(\vect{y})\Rnorm\leq M\Lnorm
\vect{x}-\vect{y}\Rnorm$, for $\vect{x},\vect{y}\in C$. For a convex
function $f$ with $M$-Lipschitz gradient, one has $ \Lnorm \nabla
f(\vect{z})-\nabla f(\vect{x})\Rnorm^2\leq M (\vect{z}-\vect{x})^\top(
\nabla f(\vect{z})-\nabla f(\vect{x}))$ for all $\vect{x},\vect{z}\in
C\subset \reals^d$.

We briefly review basic concepts from algebraic graph theory
following~\citep{FB-JC-SM:09}.  A \emph{digraph}, is a pair $\GG =
(\VV ,\EE )$, where $\VV=\{1,\dots,N\}$ is the \emph{node set} and
$\EE \subseteq \VV\times \VV$ is the \emph{edge set}.  An edge from
$i$ to $j$, denoted by $(i,j)$, means that agent $j$ can send
information to agent $i$. For an edge $(i,j) \in\EE$, $i$ is called an
\emph{in-neighbor} of $j$ and $j$ is called an \emph{out-neighbor}
of~$i$.  A graph is \emph{undirected} if $(i,j) \in \EE$ anytime
$(j,i)\in\EE$.  A \emph{directed path} is a sequence of nodes
connected by edges. A digraph is \emph{strongly connected} if for
every pair of nodes there is a directed path connecting them.  A
\emph{weighted digraph} is a triplet $\GG = (\VV ,\EE,
\vect{\sf{A}})$, where $(\VV ,\EE )$ is a digraph and $
\vect{\sf{A}}\in\real^{N\times N}$ is a weighted \emph{adjacency}
matrix such that $ \sf{a}_{ij} >0$ if $(i, j) \in\EE$ and $
\mathsf{a}_{ij} = 0$, otherwise. A weighted digraph is
\emph{undirected} if $ \mathsf{a}_{ij} = \mathsf{a}_{ji}$ for all
$i,j\in\VV$. We refer to a strongly connected and undirected graph as
a \emph{connected graph}. The \emph{weighted in-} and
\emph{out-degrees} of a node $i$ are, respectively,
$\mathsf{d}_{\text{in}}^i =\Sigma^N_{j =1} \mathsf{a}_{ji}$ and
$\mathsf{d}_{\text{out}}^i =\Sigma^N_{j =1} \mathsf{a}_{ij}$.  A
digraph is \emph{weight-balanced} if at each node $i\in\VV$, the
weighted out-degree and weighted in-degree coincide (although they
might be different across different nodes).  Any undirected graph is
weight-balanced. The \emph{(out-) Laplacian} matrix is $\lL =
\vect{\mathsf{D}}^{\text{out}} - \vect{\mathsf{A}}$, where
$\vect{\mathsf{D}}^{\text{out}} =
\Diag{\mathsf{d}_{\text{out}}^1,\cdots, \mathsf{d}_{\text{out}}^N} \in
\reals^{N \times N}$.  Note that $\lL\vect{1}_N=\vect{0}$. A digraph
is weight-balanced iff $\vect{1}_N^T\lL=\vect{0}$ iff $\Sym{\lL} =
(\vect{\lL}+\vect{\lL})/2$ is positive semi-definite. Based on the
structure of $\lL$, at least one of the eigenvalues of $\lL$ is zero
and the rest of them have nonnegative real parts.  We denote the
eigenvalues of $\lL$ by $\lambda_1, \dots,\lambda_N$, where
$\lambda_1=0$ and $\Re(\lambda_i)\leq \Re(\lambda_j)$, for $i<j$, and
the eigenvalues of $\Sym{\lL} $ by
$\hat{\lambda}_1,\dots,\hat{\lambda}_N$. For a strongly connected and
weight-balanced digraph, zero is a simple eigenvalue of both
$\vect{L}$ and $\Sym{\lL} $. In this case, we order the eigenvalues of
$\Sym{\lL} $ as $\hat{\lambda}_1 = 0 < \hat{\lambda}_2 \leq
\hat{\lambda}_3\leq\dots\leq\hat{\lambda}_N$ and we have
\vspace{-0.1in}
   \begin{equation}\label{eq::RLR}
    0<\hat{\lambda}_2 \vect{I}\leq \rR^\top\Sym{\lL}\rR\leq \hat{\lambda}_N \vect{I}. \vspace{-0.1in}   \end{equation} 
 Notice that for connected graphs $\hat{\lambda}_i=\lambda_i$ for $i\in\{1,\cdots,N\}$. For convenience, we
define $\LL = \lL\kronecker\vect{\mathsf{I}}_d$ and $\PPi=
\pPi_N\kronecker\vect{\mathsf{I}}_d$ to deal with variables of
dimension $d\in \naturals$.

\section{Problem Definition}\label{se:problem-statement}

Consider a network of $N$ agents interacting over~a strongly connected and
weight-balanced digraph~$\GG$. Each agent $i
\until{N}$~has a differentiable local cost function $f^i:\mathbb{R}^d\to\mathbb{R}$. 
The global cost function of the network, which we assume to be strictly convex, is $f=\Sigma_{i=1}^N f^i(\vect{\mathsf{x}})$. 
Our aim is to design a distributed algorithm such that each
agent solves
\vspace{-0.05in}
\begin{equation*}
  \min_{\vect{\mathsf{x}}\in
    \reals^d} \,\,f(\vect{\mathsf{x}}) 
\end{equation*}
using only its own local data and exchanged information with its
neighbors. We assume the above optimization problem is feasible
(which, together with the strict convexity of $f$, implies the
uniqueness of the global optimizer, which we denote
$\vect{\mathsf{x}}^{\star} \in \reals^d$).  We are also interested in
characterizing the privacy preservation properties of the algorithmic
solution to this distributed optimization problem. Specifically, we
aim to identify conditions guaranteeing that no information about the
local cost function of an agent is revealed to, or can be
reconstructed by, any other agent in the network.

\section{Distributed Continuous-Time Algorithm for Convex
  Optimization}\label{sec::Alg}

In this section, we provide a novel continuous-time distributed
coordination algorithm to solve the problem stated in
Section~\ref{se:problem-statement} and analyze in detail its
convergence properties. For $i\!\until{N}$ and with $\alpha,\beta\!>\!0$, consider
\vspace{-0.1in}
\begin{subequations}\label{eq::Alg}
  \begin{align}
    &\dvect{v}^i = \alpha\beta\sum\nolimits_{j=1}^N
    \mathsf{a}_{ij}(\vect{x}^i-\vect{x}^j), \label{eq::Alg-b}
    \\
    &\dvect{x}^i = -\alpha\nabla f^i(\vect{x}^i) - \beta\sum\nolimits_{j=1}^N
    \mathsf{a}_{ij}(\vect{x}^i-\vect{x}^j) - \vect{v}^i, \label{eq::Alg-a}
  \end{align}
\end{subequations}
In network variables $\vect{x},\!\vect{v} \!\in\! (\real^d)^N\!$, the algorithm reads as
\vspace{-0.1in}
\begin{subequations}\label{eq::Alg_col}
  \begin{align}
    \dvect{v} & = \alpha\beta\LL\vect{x},\label{eq::Alg_col_v}
    \\
    \dvect{x} & =-\alpha\nabla
    \tilde{f}(\vect{x})-\beta\LL\vect{x}-\vect{v}.\label{eq::Alg_col_x}
  \end{align}
\end{subequations}
Here, $\map{\tilde{f}}{(\real^d)^N}{\real}$ is defined by $
\tilde{f}(\vect{x})=\Sigma_{i=1}^Nf^i(\vect{x}^i)$.  This algorithm is
distributed because each agent only needs to receive information from its
out-neighbors about their corresponding variables in $\vect{x}$. In
contrast, the continuous-time coordination algorithms
in~\citep{JW-NE:11,BG-JC:14-tac} require the communication of the
corresponding variables in both $\vect{x}$ and~$\vect{v}$.
The synthesis of algorithm~\eqref{eq::Alg_col} is inspired by the
following feedback control considerations.  In~\eqref{eq::Alg_col_x},
each agent follows a local gradient descent while trying to agree with
its neighbors on their estimate of the final value. However, as the
local gradients are not the same, this dynamics by itself would never
converge. Therefore, to correct this error, each agent uses an
integral feedback term $\vect{v}^i$ whose evolution is driven by the
agent disagreement according to~\eqref{eq::Alg_col_v}.  

Our analysis of the algorithm convergence is structured in two parts,
depending on the directed character of the
interactions. Section~\ref{sec:directed} deals with strongly connected
and weight-balanced digraphs and Section~\ref{sec:undirected} deals
with connected graphs. In each case, we identify conditions
on the agent cost functions that guarantee asymptotic
convergence. Given the challenges posed by directed information flows,
it is not surprising that we can establish stronger results under less
restrictive assumptions for the case of undirected topologies.

\subsection{Strongly Connected, Weight-Balanced
  Digraphs}\label{sec:directed}

Here, we study the convergence of the distributed optimization
algorithm~\eqref{eq::Alg} over strongly connected and weight-balanced
digraph topologies. We first consider the case where the interaction
topology is fixed, and then discuss the time-varying interaction
topologies.  The following result identifies conditions on the local
cost functions $\{f^i\}_{i=1}^N$ and the parameter $\beta$ to
guarantee the exponential convergence of~\eqref{eq::Alg} to the
solution of the distributed optimization problem.

\begin{thm}[Convergence of~\eqref{eq::Alg} over strongly connected and
  weight-balanced digraphs]\label{thm::OpSlvrDirected}
  Let $\GG$ be a strongly connected and weight-balanced digraph.
  Assume each $f^i$, $i\!\until{N}$, is $m^i$-strongly convex,
  differentiable, and its gradient is $M^i$-Lipschitz
  on~$\real^d$. Given $\alpha>0$, $\mT = \min\{m^1,\dots,m^N\}$ and $
  \MT =\max\{M^1,\dots,M^N\}$, let $\beta,\phi>0$ satisfy $ \phi+1>
  4\MT$ and \vspace{-0.1in}
  \begin{align}\label{eq::gamma}
    \gamma\!= \alpha^2(\phi\!+\!1)\mT\! +\!
    9\beta\Hlambda_2\phi\alpha \!-\!  4\alpha^2(
    \MT\mT\!+\!(\phi\!+\!1)^2) \!\!>\!0, 
  \end{align}
  Then, for each $i\until{N}$, starting from 
  $\vect{x}^i(0),\vect{v}^i(0)\in\reals^d$ with
  $\sum_{i=1}^N\vect{v}^i(0) \!=\! \vect{0}_d$, the
  algorithm~\eqref{eq::Alg} over $\GG$ makes
  $\vect{x}^i(t)\to\vect{\mathsf{x}}^\star$ exponentially fast as
  $t\to\infty$ with a rate no less than 
    \begin{equation}\label{eq::rate_expon}
    \min\{ \frac{7}{16},\frac{1}{9}\gamma\}/(2\bar{\lambda}_{\vect{F}}).
 \end{equation}
 Here, $\bar{\lambda}_{\vect{F}} $ is the maximum eigenvalue of
 \vspace{-0.1in}
  \begin{equation}\label{eq:matrix-F}
    \vect{F}=\frac{1}{2}\!\!
    \begin{bmatrix}
      \frac{1}{9}\alpha(\phi+1)\vect{\mathsf{I}}_d & \vect{0} &
      \vect{0}
      \\
      \vect{0}& \alpha(\phi+1)\vect{\mathsf{I}}_{(N-1)d} &
      \vect{\mathsf{I}}_{(N-1)d}
      \\
      \vect{0}&\vect{\mathsf{I}}_{(N-1)d} &
      \frac{1}{\alpha}\vect{\mathsf{I}}_{(N-1)d}
    \end{bmatrix}.
  \end{equation}
\end{thm}
\begin{pf}
  For weight-balanced digraphs, we have
  $\vect{1}_N^\top\lL\!=\!\vect{0}$. Thus,
  left multiplying~\eqref{eq::Alg_col_v} by $\vect{1}_N^\top\!\kronecker\!
  \vect{\mathsf{I}}_d$ results~in
  \vspace{-0.1in}
  \begin{align}\label{eq::sum_v}
   \! \!\sumSig{\dot{\vect{v}}^i} \!=\! \vect{0} \Rightarrow\!
    \sumSig{\vect{v}^i(t)} \!=\! \sumSig{\vect{v}^i(0)} \!=\!
    \vect{0},~\forall t\geq 0.
  \end{align}
  Next, we obtain the equilibrium point of \eqref{eq::Alg},
  ($\Bvect{x}, \Bvect{v}$), from 
  \vspace{-0.1in}
  \begin{subequations}\label{eq::Alg_col_Equilib}
    \begin{align}
      &\vect{0} =
      \alpha\beta\LL\Bvect{x},\label{eq::Alg_col_Equilib-b}
      \\
      &\vect{0} = -\alpha \nabla
      \tilde{f}(\Bvect{x})-\beta\LL\Bvect{x}-\Bvect{v}.
      \label{eq::Alg_col_Equilib-a}
    \end{align}
  \end{subequations}
  Given~\eqref{eq::Alg_col_Equilib-b}, $\Bvect{x}$
  belongs to the null-space of $\LL $. For strongly connected digraphs 
  the null-space of $\lL$ is 
  spanned by
  $\vect{1}_{N}$. Thus, $\Bvect{x} = \vect{1}_{N}\kronecker
    \vect{\theta}$, where $\vect{\theta}\in\reals^d$.
  Left multiplying~\eqref{eq::Alg_col_Equilib-a} by
  $\vect{1}_N^\top\kronecker \vect{\mathsf{I}}_d$ and
  using \eqref{eq::sum_v}, we obtain
  $\vect{0}=\sum_{i=1}^{N} \nabla f^i(\bar{\vect{x}}^i)$. 
  Then, the optimality condition $ \nabla f(\vect{\mathsf{x}}^\star) =
  \vect{0}_d$ along with $\nabla
  f(\vect{\theta}) =\sum_{i=1}^N \nabla f^i(\vect{\theta})$ imply
  \vspace{-0.05in}
  \begin{equation*}\label{eq::x_bar}
    \bar{\vect{x}}^i=\vect{\mathsf{x}}^\star,~~~i\until{N}.
  \end{equation*}
  Substituting this value in \eqref{eq::Alg_col_Equilib-a}, we obtain
  \vspace{-0.1in}
  \begin{equation}\label{eq::v_bar}
    \bar{\vect{v}}^i=-\alpha\nabla f^i(\vect{\mathsf{x}}^\star),~~~i\until{N}.
  \end{equation}
  To study the stability of~\eqref{eq::Alg}, we transfer
  the equilibrium point to the origin and then apply a change of variables
  \vspace{-0.1in}
  \begin{subequations}
  \begin{align}
   & \vect{u}=\vect{v}-\Bvect{v},\quad\vect{y}=\vect{x}-\Bvect{x},\label{eq::shift_orig}\\
   & \vect{u}= ([
     \vect{\mathsf{r}}\quad
     \rR
     ]
     \kronecker
     \vect{\mathsf{I}}_d)\vect{w},
    \quad\vect{y}=([
     \vect{\mathsf{r}}\quad
     \rR
     ]
     \kronecker
     \vect{\mathsf{I}}_d)\vect{z},\label{eq::trans_stable}
  \end{align}
  \end{subequations}
  where we used~\eqref{eq::Rrs}. We partition the new variables~as follows:
  $\vect{w}=(\vect{w}_1, \vect{w}_{2:N})$ and $\vect{z}=(\vect{z}_1,
  \vect{z}_{2:N})$, where $\vect{w}_1,\vect{z}_1\in\reals^d$.
  In these new variables, the algorithm~\eqref{eq::Alg} reads as
  \vspace{-0.1in}
  \begin{align}\label{eq::Alg_Separated}
    \dvect{w}_{1} & = \vect{0}_d, \notag
    \\
    \dvect{w}_{2:N} & = \alpha\beta(\rR^\top\lL\rR \kronecker
    \vect{\mathsf{I}}_d)\vect{z}_{2:N}, 
    \notag
    \\
    \dvect{z}_{1} & =- \alpha(\vect{\mathsf{r}}^\top\kronecker
    \vect{\mathsf{I}}_d)\vect{h}, 
    \\
    \dvect{z}_{2:N} & = -\alpha(\rR^\top\kronecker
    \vect{\mathsf{I}}_d)\vect{h}-\beta(\rR^\top\lL\rR \kronecker
    \vect{\mathsf{I}}_d)\vect{z}_{2:N}-\vect{w}_{2:N}
    , 
    \notag
  \end{align}
  where 
  \vspace{-0.1in}
  \begin{align}\label{eq::h}
    \vect{h} = \nabla \tilde{f}(\vect{y}+\Bvect{x})-\nabla
    \tilde{f}(\Bvect{x}).
  \end{align}
  Note that the first equation in~\eqref{eq::Alg_Separated}
  corresponds to the constant of motion~\eqref{eq::sum_v}.  To study
  the stability in the other variables, consider the candidate
  Lyapunov function
  \vspace{-0.1in}
  \begin{align}\label{eq::Lya_Digraph}
    V(\vect{z},\vect{w}_{2:N}) &=
    \frac{1}{18}\alpha(\phi+1){\vect{z}_1}^\top\vect{z}_1 +
    \frac{\phi\alpha}{2}{\vect{z}_{2:N}}^\top\vect{z}_{2:N}
    \\
    & \; + \frac{1}{2\alpha}(\alpha\vect{z}_{2:N}+\vect{w}_{2:N})^\top
    (\alpha\vect{z}_{2:N}+\vect{w}_{2:N}) , \nonumber
  \end{align}
  with $\phi\!>\!0$ as in the statement. Note that $
  V \!\leq \!\bar{\lambda}_{\vect{F}}
  \|\vect{p}\|^2$, with $\vect{p}\!=\!(\vect{z},\vect{w}_{2:N})$.
  The Lie derivative of $V$ along~\eqref{eq::Alg_Separated} is 
  \vspace{-0.1in}
  \begin{align*}
    \!\! \dot{V} =& -\frac{1}{9}\alpha^2(\phi+1)\vect{y}^\top\vect{h}
    - \frac{7}{16}\vect{w}_{2:N}^\top\vect{w}_{2:N}
    \\
    &-\phi\alpha\beta{\vect{z}_{2:N}}^\top(\rR^\top
    \Sym{\lL} \rR \kronecker \vect{\mathsf{I}}_d)\vect{z}_{2:N}
    \\
    &+\frac{4}{9} \alpha^2\|(\rR^\top\kronecker
    \vect{\mathsf{I}}_d)\vect{h}\|^2 + \frac{4}{9}
    \alpha^2(1+\phi)^2\vect{z}_{2:N}^\top\vect{z}_{2:N}
    \\
    &-\|\frac{3}{4}\vect{w}_{2:N}+\frac{2\alpha}{3}(\rR^\top\kronecker
    \vect{\mathsf{I}}_d)\vect{h}+\frac{2\alpha}{3}(\phi+1)\vect{z}_{2:N})\|^2.
  \end{align*} 
  Next, we show that under $ \phi+1> 4\MT$
  and~\eqref{eq::gamma}, $\dot{V}$ is negative
  definite. Invoking the
  assumptions on the local cost functions in the statement, and using
  the $\MT$-Lipschitzness of $\nabla \tilde{f}$ and the $\mT$-strongly
  convexity of $\tilde{f}$ along with $\|\rR^\top\kronecker
  \vect{\mathsf{I}}_d \|=1$ and $ \|\vect{z}\|=\|\vect{y}\|$, we have \vspace{-0.1in}
  \begin{subequations}\label{eq::mT-MT}
    \begin{align}
      & \|(\rR^\top\kronecker \vect{\mathsf{I}}_d)\vect{h}\|^2 \leq
      \|\vect{h}\|^2\leq \MT\vect{y}^\top\vect{h},\label{eq::mT-MT-M}
      \\
      & {\vect{y}}^\top\vect{h} \geq \mT \|\vect{y}\|^2 =\mT
      \|\vect{z}\|^2.\label{eq::mT-MT-m}
    \end{align}
  \end{subequations} 
  Given these
  relations and invoking~\eqref{eq::RLR}, we have \vspace{-0.1in}
  \begin{align*}
    \dot{V} \!& \leq \!-\frac{\alpha^2((\phi\!+\!1)- 4 \MT)
      \mT}{9}\vect{z}^\top\vect{z} \!-\!
    \frac{7}{16}\vect{w}_{2:N}^\top\vect{w}_{2:N}\! 
    \\
    & \quad -
    \phi\alpha\beta\Hlambda_2\vect{z}_{2:N}^\top\vect{z}_{2:N}+
    \frac{4\alpha^2(1+\phi)^2}{9}\vect{z}_{2:N}^\top\vect{z}_{2:N}\\
    & \quad
    -\|\frac{3}{4}\vect{w}_{2:N}+\frac{2\alpha}{3}(\rR^\top\kronecker
    \vect{\mathsf{I}}_d)\vect{h}+\frac{2\alpha(\phi+1)}{3}\vect{z}_{2:N})\|^2.
  \end{align*}
  Since $\vect{z}^\top\vect{z} = \vect{z}_1^\top\vect{z}_1 +
  \vect{z}_{2:N}^\top\vect{z}_{2:N}$, it follows that $\dot{V}<-\min\{
  \frac{7}{16},\frac{1}{9}\gamma\}\|\vect{p}\|^2<0$, where $\gamma$ is
  a shorthand notation for the expression
  in~\eqref{eq::gamma}.  Thus, 
  $\vect{z}\to\vect{0}$ as $t\to\infty$, equivalently
  $\vect{x}^i\to\vect{\mathsf{x}}^\star$, for all $i\until{N}$, is
  exponential with rate no less than~\eqref{eq::rate_expon}
  (cf.~\cite[Theorem 4.10]{HKK:02}).  \boxend
\end{pf}

In Theorem~\ref{thm::OpSlvrDirected}, the requirement $\sum_{i=1}^N
\vect{v^i}(0)=\vect{0}_d$ is trivially satisfied by each agent with
$\vect{v}^i(0)=\vect{0}_d$. This is an advantage with respect to the
continuous-time coordination algorithms proposed in~\citep{JL-CYT:12},
which requires the nontrivial initialization $\sum_{i=1}^{N} \nabla
f^i(\vect{x}^i(0)) = \vect{0}_d$, and in~\citep{FZ-DV-AC-GP-LS:11},
which requires the initialization on a state communicated among
neighbors and is subject to communication error.

\begin{rem}\longthmtitle{Design parameters
    in~\eqref{eq::Alg}}\label{rem:design-parameters} {\rm We provide
    here several observations on the role of the design
    parameters $\alpha$ and $\beta$. First, note  there always
    exist $\alpha$, $\beta$
    satisfying~\eqref{eq::gamma}, e.g., any $\beta>
    {4(\phi+1)^2\alpha}/{(9\phi\Hlambda_2)}$.
    The determination of these parameters can be performed by
    individual agents if they know an upper bound on $\MT$, a lower
    bound on $\mT$, and have knowledge of $\Hlambda_2$, either through
    a dedicated algorithm to compute it, see
    e.g.,~\cite{PY-RAF-GJG-KML-SSS-RS:10}, or use a lower bound on it,
    see e.g.~\cite{BM:91a}.  We have observed in simulation
    that~\eqref{eq::gamma} is only sufficient and that,
    in fact, the algorithm~\eqref{eq::Alg} converges for any positive
    $\alpha$ and $\beta$ in our numerical examples.  Although not
    evident in~\eqref{eq::rate_expon}, one can expect the larger
    $\alpha$ and $\beta$ are, the higher the rate of convergence of
    the algorithm~\eqref{eq::Alg} is. A coefficient $\alpha>1$ can be
    interpreted as a way of increasing the strong convexity
    coefficient of the local cost functions. A coefficient $\beta>1$
    can be interpreted as a means of increasing the graph
    connectivity.  The relationship between these parameters and the
    rate of convergence of~\eqref{eq::Alg} is more evident for
    quadratic local cost functions $f^i(\vect{\mathsf{x}}) =
    \frac{1}{2}(\vect{\mathsf{x}}^\top\vect{\mathsf{x}} +
    \vect{\mathsf{x}}^\top \vect{a}^i+\vect{b}^i)$, $i\until{N}$. In
    this case, the algorithm~\eqref{eq::Alg} is a linear
    time-invariant system where the eigenvalues of the system matrix
    are $-\alpha$, with multiplicity of $Nd$, and $-\beta\lambda_i$,
    $i\until{N}$ ($\lambda_i$'s are the eigenvalues of $\lL$), with
    multiplicity $d$.  Therefore,~\eqref{eq::Alg} converges regardless
    of the value of $\alpha$, $\beta>0$ with an exponential rate equal
    to $\min\{\alpha,\beta\Re(\lambda_2)\}$.  When discussing
    discrete-time communication, some trade-offs arise regarding the
    choice of the parameters, as we explain later in
    Section~\ref{sec::Discrete}.}  \oprocend
\end{rem}

\begin{rem}\longthmtitle{Semiglobal convergence of~\eqref{eq::Alg} under local
  gradients that are locally Lipschitz}\label{rem::Lips_compact}
  {\rm The convergence result in Theorem~\ref{thm::OpSlvrDirected} is
    semiglobal~\citep{HKK:02} if the local gradients are only locally
    Lipschitz or, equivalently, Lipschitz on compact sets. In fact,
    one can see from the proof of the result that, for any compact set
    containing the initial conditions $\vect{x}^i(0)\in\real^d$ and
    $\vect{v}^i(0)=\vect{0}_d$, $i\until{N}$, one can find $\phi>0$
    and $\beta>0$ sufficiently large such that the compact set is
    contained in the region of attraction of the equilibrium
    point. \oprocend }
\end{rem}
 
Next, we study the convergence of~\eqref{eq::Alg} over dynamically
changing topologies. Since the proof of
Theorem~\ref{thm::OpSlvrDirected} relies on a Lyapunov function with
no dependency on the system parameters and its derivative is upper
bounded by a quadratic negative definite function, we can readily
extend the convergence result to dynamically changing networks. The
proof details are omitted for brevity.

\begin{prop}[Convergence of~\eqref{eq::Alg} over dynamically changing
  interaction topologies]\label{prop::OpSlvrDirectedSwtch}
  Let $\GG$ be a time-varying digraph which is strongly connected and
  weight-balanced at all times and whose adjacency matrix is uniformly
  bounded and piecewise constant. Assume the local cost function
  $f^i$, $i\!\until{N}$, is $m^i$-strongly convex, differentiable, and
  its gradient is $M^i$-Lipschitz on~$\real^d$.  Given $\alpha>0$, let
  $\beta,\phi>0$ satisfy $ \phi+1> 4\MT$
  and~\eqref{eq::gamma} with $\Hlambda_2$ replaced by
  $(\Hlambda_2)_{\min} \!=\!
  \underset{p\in\PP}{\min}\{\Hlambda_2(\lL_p)\}$, where $\PP$ is the
  index set of all possible realizations of $\GG$.  Then, for each
  $i\until{N}$, starting from 
  $\vect{x}^i(0),\vect{v}^i(0)\!\in\!\reals^d$ with
  $\sum_{i=1}^N\!\vect{v}^i(0)\!=\!\vect{0}_d$, the
  algorithm~\eqref{eq::Alg} over $\GG$ makes
  $\vect{x}^i(t)\to\vect{\mathsf{x}}^\star$ exponentially fast as
  $t\to\infty$.
\end{prop}

Our final result of this section characterizes the topological
requirements on the communication graph and the knowledge about the
algorithm's parameters and initial conditions that allow a passive
agent (i.e., an agent that does not interfere in the algorithm
execution) to reconstruct the local gradients of other agents in the
network.

\begin{prop}[Privacy preservation
  under~\eqref{eq::Alg}]\label{prop::privLocalg}
  Let $\GG$ be a strongly connected and weight balanced digraph. For
  $\alpha,\beta>0$, consider any execution of the coordination
  algorithm~\eqref{eq::Alg} over $\GG$ starting from
  $\vect{x}^i(0),\vect{v}^i(0)\in\real^d$ with
  $\sum_{i=1}^{N}\vect{v}^i(0) = \vect{0}_d$. Then, an agent $i
  \until{N}$ can reconstruct the local gradient of another agent $j
  \neq i$ only if $j$ and all its out-neighbors are out-neighbors of
  $i$, and agent $i$ knows $\vect{v}^j(0)$ and $ a_{jk}$, $k\until{N}$
  (here we assume that the agent $i$ is aware of the identity of
  neighbors of agent $j$ and it has memory to save the time history of
  the data it receives from its out-neighbors).
\end{prop}
\begin{pf} 
  Consider an arbitrary time $t^\star$.  Let $i$ be an in-neighbor of
  agent $j$ and all of its out-neighbors. The
  algorithm~\eqref{eq::Alg} requires each agent to communicate its
  component of $\vect{x}$ to their in-neighbors. Since agent $i$ has
  memory to save information it receives from its out-neighbors for
  all $t\leq t^\star$, it can use the time history of $\vect{x}^j(t)$
  to numerically reconstruct $\dvect{x}^j(t^\star)$. Because $i$ is
  the in-neighbor of $j$ and its out-neighbors, it can use its
  knowledge of $a_{jk}$, $k\until{N}$ to reconstruct $\sum_{k=1}^N\!
  a_{jk}(\vect{x}^j(t) - \vect{x}^k(t))$ for all $t\leq
  t^\star$. Agent $i$ can reconstruct $\vect{v}^j(t)$
  from~\eqref{eq::Alg-b} uniquely as it knows $\vect{v}^j(0)$. Then,
  agent $i$ has all the elements to solve for $\nabla
  f^j(\vect{x}^j(t^\star))$ in~\eqref{eq::Alg-a}.  The lack of
  knowledge about any of this information would prevent $i$ from
  reconstructing exactly the local gradient of~$j$.  \boxend
\end{pf}

\subsection{Connected Graphs}\label{sec:undirected}

Here, we study the convergence of the algorithm~\eqref{eq::Alg} over
connected graph topologies. While the results of the previous section
are of course valid for these topologies, here, using the structural
properties of the Laplacian matrix of undirected graphs, we establish
the convergence of~\eqref{eq::Alg} for a larger family of local cost
functions. In doing so, we are also able to analytically establish
convergence for any $\alpha,\beta>0$, as we show next.
\begin{thm}[Exponential convergence of~\eqref{eq::Alg} over connected
  graphs]\label{thm::UD_ConvRate_Alg}
  Let $\GG$ be a connected graph. Assume the local cost function
  $f^i$, $i\!\until{N}$, is $m^i$-strongly convex and differentiable
  on $\real^d$, and its gradient is locally Lipschitz.  Then, for any
  $\alpha,\beta>0$ and each $i\until{N}$, starting from  
  $\vect{x}^i(0),\vect{v}^i(0)\in\reals^d$ with
  $\sum_{i=1}^N\vect{v}^i(0) = \vect{0}$, the
  algorithm~\eqref{eq::Alg} over $\GG$ satisfies
  $\vect{x}^i(t)\to\vect{\mathsf{x}}^\star$ as $t\to\infty$
  exponentially fast.
\end{thm}
\begin{pf}
  We use the equivalent representation~\eqref{eq::Alg_Separated} of
  the algorithm~\eqref{eq::Alg} obtained in the proof of
  Theorem~\ref{thm::OpSlvrDirected}.  
  Consider the following candidate
  Lyapunov function
  \vspace{-0.1in}
  \begin{align}\label{eq::LyapunovStronConvexUndi}
    V(\vect{z},\vect{w}_{2:N}) & =
    \frac{1}{2}\alpha(\phi+1)\vect{z}_1^\top\vect{z}_1 +
    \frac{\phi\alpha}{2} \vect{z}_{2:N}^\top\vect{z}_{2:N}
    \\
    & \quad +
    \frac{1}{2\alpha}(\alpha\vect{z}_{2:N}+\vect{w}_{2:N})^\top
    (\alpha\vect{z}_{2:N}+\vect{w}_{2:N}) \nonumber
    \\
    & \quad + \frac{1}{2\beta}(\phi+1)\vect{w}_{2:N}^\top
    \big((\rR^\top\lL\rR )^{-1} \kronecker
    \vect{\mathsf{I}}_d\big) \vect{w}_{2:N}, \nonumber
  \end{align}
  with $\phi\!\geq\! 1$ defined below. Given~\eqref{eq::RLR}, $V$ is
  positive definite and radially unbounded and satisfies $V\! \leq\!
  \bar{\lambda}_{\vect{E}} \|\vect{p}\|^2$, with
  $\vect{p}\!=\!(\vect{z},\vect{w}_{2:N})$ and
  $\bar{\lambda}_{\vect{E}}\!>\!0 $ is the maximum eigenvalue of
  \vspace{-0.1in} {\setlength\arraycolsep{1pt}
    \begin{equation*}
      \vect{E}\!=\!\frac{1}{2}\!\!
      \begin{bmatrix}
        \!\alpha(\phi\!+\!1)\vect{\mathsf{I}}_d\!&\!\vect{0}&\vect{0}
        \\
        \vect{0}&
        \!\alpha(\phi\!+\!1)\vect{\mathsf{I}}_{(N-1)d}\!\!&\vect{\mathsf{I}}_{(N-1)d}
        \\
        \vect{0}&\vect{\mathsf{I}}_{(N-1)d}&\frac{1}{\alpha}\vect{\mathsf{I}}
        \!+\!\frac{(\phi\!+\!1)}{\beta}(\rR^\top\!\lL\rR )^{-1}\! \kronecker\!
        \vect{\mathsf{I}}_d\!
      \end{bmatrix}\!\!.
    \end{equation*}}
  The Lie derivative of $V$ along the
  dynamics~\eqref{eq::Alg_Separated} is
  \vspace{-0.1in}
   \begin{align*}
    \dot{V} & =-\alpha^2(\phi+1) \vect{y}^\top\vect{h}- \phi \alpha
    \beta{\vect{z}_{2:N}}^\top(\rR^\top\lL\rR \kronecker
    \vect{\mathsf{I}}_d)\vect{z}_{2:N}
    \\
    & \quad - {\vect{w}_{2:N}}^\top
    \vect{w}_{2:N}-\alpha{\vect{w}_{2:N}}^\top(\rR^\top\kronecker
    \vect{\mathsf{I}}_d)\vect{h} . 
  \end{align*}
  Next, we show that $\dot{V}$ is upper bounded by a
  negative definite function.  We start by identifying a
  compact set whose definition is independent of $\phi$ and contains
  the set $\SLya=\{ (\vect{z},\vect{w}_{2:N}) \in \real^{Nd} \times
  \real^{(N-1)d}   \; | \;
  V(\vect{z},\vect{w}_{2:N}) \leq
  V(\vect{z}(0),\vect{w}_{2:N}(0))\}$.    
  For any given initial condition, let $\rho_0=
  \frac{1}{2}\alpha\vect{z}_1(0)^\top\vect{z}_1(0) +
  \frac{(\alpha+1)}{2} \vect{z}_{2:N}(0)^\top\vect{z}_{2:N}(0)+ (
  \frac{1}{2\beta\lambda_2}+\frac{1}{2\alpha}+\frac{1}{2})
  \vect{w}_{2:N}(0)^\top\vect{w}_{2:N}(0)$ and define $\SBar
   =\{ \vect{z} \in \real^{Nd} \; | \;
  \frac{1}{2}\alpha\vect{z}_1^\top\vect{z}_1 +
  \frac{1}{4}\alpha\vect{z}_{2:N}^\top\vect{z}_{2:N}\leq \rho_0\}$.
  Observe that this set is compact.  
  Note that $(\vect{z},\vect{w}_{2:N}) \in \SLya$ implies
  $\vect{z}\in\SBar$ because $ V(\vect{z}(0),\vect{w}_{2:N}(0)) \leq
  (\phi+1)\rho_0$ and
  \begin{align*}
    & \frac{1}{2}\alpha(\phi+1) \vect{z}_1^\top\vect{z}_1 +
    \frac{1}{4}\alpha(\phi+1)\vect{z}_{2:N}^\top\vect{z}_{2:N}\leq
    \\
    & \frac{1}{2}\alpha(\phi+1)\vect{z}_1^\top\vect{z}_1 +
    \frac{1}{2}\alpha\phi \vect{z}_{2:N}^\top\vect{z}_{2:N}\leq
    V(\vect{z},\vect{w}_{2:N}).
  \end{align*}
  Here, we used $\phi\geq1$ in the first inequality.  Since the
  change of variables~\eqref{eq::shift_orig}
  and~\eqref{eq::trans_stable} are linear, the corresponding
  $\vect{x}$ and $\vect{y}$ for $\vect{z}\in\SBar$ belong to
  compact sets, as well.  Then, the assumption on the gradients of the
  local cost functions implies that there exists $M_0>0$ with   
   $ \|\vect{h}\|\leq M_0 \|\vect{y}\|=M_0 \|\vect{z}\|$, for all $
    \vect{z}\in \SBar$. 
  Consequently,
  \begin{equation*}
    \|(\rR^\top\kronecker \vect{\mathsf{I}}_d) \vect{h}\|^2\!\leq M_0^2
    \|\vect{z}\|^2,\quad 
    \forall (\vect{z},\vect{w}_{2:N})  \in \SLya. 
  \end{equation*}
Then we can show
  $-\alpha{\vect{w}_{2:N}}^\top(\rR^\top\kronecker
  \vect{\mathsf{I}}_d)\vect{h}\leq \frac{1}{2}\vect{w}_{2:N}^\top\vect{w}_{2:N}
  + \frac{1}{2}\alpha^2M_0^2\vect{z}^\top\vect{z}$ for
  $(\vect{z},\vect{w}_{2:N}) \in \SLya$. Using this inequality, 
  fact that the local cost functions are $m^i$-strongly convex (and
  hence~\eqref{eq::mT-MT-m} holds) and invoking~\eqref{eq::RLR} 
  we
  deduce
  \vspace{-0.1in}
  \begin{align*}
    \dot{V}& \leq -\alpha^2(\phi+1)\mT(\vect{z}_1^\top\vect{z}_1 +
    \vect{z}_{2:N}^\top\vect{z}_{2:N})-\phi\alpha\beta \lambda_2
    \vect{z}_{2:N}^\top\vect{z}_{2:N} \nonumber
    \\
    & \quad - \frac{1}{2} \vect{w}_{2:N}^\top\vect{w}_{2:N}+
    \frac{1}{2}\alpha^2 M_0^2(\vect{z}_1^\top\vect{z}_1 +
    \vect{z}_{2:N}^\top\vect{z}_{2:N}).
  \end{align*}
  Let $ \phi+1\! =\!
  \frac{1}{2\mT}M_0^2\!+\!\frac{1}{2\mT\alpha^2}\delta_0$, where
  $\delta_0\!>\!0$ is such that $\phi\!\geq \!1$ (since $M_0$ does not
  depend on $\phi$, this choice is always feasible). Then, we have
  $\dot{V}\!\leq\! - \frac{1}{2}\min\{1,\delta_0\}\|\vect{p}\|^2$.~As
  such, the Lyapunov function~\eqref{eq::LyapunovStronConvexUndi}
  satisfies all the conditions of~\cite[Theorem~4.10]{HKK:02}, for the
  dynamics~\eqref{eq::Alg_Separated}.  Therefore, the convergence of
  $\vect{z}\to\vect{0}$ as $t\to\infty$, equivalently
  $\vect{x}^i\to\vect{\mathsf{x}}^\star$, for all $i\until{N}$, is
  exponential. \boxend
\end{pf}

\begin{rem}\longthmtitle{Bound on exponential rate of
    convergence}\label{rem::rate_cov_contiu} 
  {\rm One can see from the proof of
    Theorem~\ref{thm::UD_ConvRate_Alg} that, for a given initial
    condition, the rate of convergence is at least
    $\frac{1}{4}(\min\{1,\delta_0\})/\bar{\lambda}_{\vect{E}}>0$. The guaranteed rate of convergence is
    therefore not uniform, unless the local gradients are globally
    Lipschitz. In this case, one recovers the result in
    Theorem~\ref{thm::OpSlvrDirected} but for arbitrary~$\alpha$,
    $\beta>0$.\oprocend}
\end{rem}

Note that in Theorem~\ref{thm::UD_ConvRate_Alg} the local
Lipschitzness of $\nabla f^i $ is trivially held if $f^i$ is twice
differentiable. The Lyapunov
function~\eqref{eq::LyapunovStronConvexUndi} identified in the proof
of this result plays a key role later in our study of the algorithm
implementation with discrete-time communication in
Section~\ref{sec::Discrete}.  Next, we study the convergence of the
algorithm~\eqref{eq::Alg} over connected graphs when the local cost
functions are only convex.  Here, the lack of strong convexity makes
us rely on a LaSalle function, rather than on a Lyapunov one, to
establish asymptotic convergence to the optimizer.

\begin{thm}[Asymptotic convergence of~\eqref{eq::Alg} over connected
  graphs]\label{thm::UD_Conv_Alg2}
  Let $\GG$ be a connected graph.  Assume the local cost function
  $f^i$, $i\!\until{N}$, is convex and differentiable on $\real^d$,
  and the global cost function $f$ is strictly convex.  Then, for any
  $\alpha,\beta>0$ and each $i\until{N}$, starting from 
$\vect{x}^i(0),\vect{v}^i(0)\in \reals^d$ with
  $\sum_{i=1}^N\vect{v}^i(0) = \vect{0}_d$, the
  algorithm~\eqref{eq::Alg} over $\GG$ make
  $\vect{x}^i(t)\to\vect{\mathsf{x}}^\star$ as $t\to\infty$.
\end{thm}
\begin{pf}
  We use again the equivalent representation~\eqref{eq::Alg_Separated}
  of the algorithm~\eqref{eq::Alg} obtained in the proof of
  Theorem~\ref{thm::OpSlvrDirected}.  To study the stability of this
  system, consider the following candidate Lyapunov function
  \vspace{-0.1in}
  \begin{align*}
        V(\vect{z},\vect{w}_{2:N}) =& \frac{1}{2}\vect{z}^\top\vect{z}
   \! +\! 
    \frac{1}{2\alpha\beta}{\vect{w}_{2:N}}^\top\!
    \left(\!(\rR^\top\lL\rR )^{-1} \!\kronecker\!
      \vect{\mathsf{I}}_d\right)\!\vect{w}_{2:N}.
  \end{align*} 
  Given~\eqref{eq::RLR}, $V$ is
  positive definite and radially unbounded. The Lie derivative of $V$
  along~\eqref{eq::Alg_Separated} is given by
  \vspace{-0.1in}
   \begin{align*}
    \dot{V} & = -\alpha\vect{y}^\top(\nabla_\text{T}
    f(\vect{y}+\Bvect{x})-\nabla_\text{T}f(\Bvect{x}))
    -\beta\vect{y}^\top\LL \vect{y}
    \\
    & = -\alpha\sum\nolimits_{i=1}^N{\vect{y}^i}^\top(\nabla f^i(\vect{y}^i
    \!+\! \vect{\mathsf{x}}^\star)-\nabla f^i(\vect{\mathsf{x}}^\star)) -
    \beta\vect{y}^\top\LL \vect{y}.
  \end{align*}
  To obtain the second summand, we have used
  $\vect{z}_{2:N}=(\rR^\top\kronecker \vect{\mathsf{I}}_d)\vect{y}$
  and $\rR \rR^\top = \vect{\mathsf{I}}_N -
  \vect{\mathsf{r}}\vect{\mathsf{r}}^\top$. Since the local cost
  functions are convex, the first summand of $\dot{V}$ is non-positive
  for all $\vect{y}$. Because the graph is connected, the second
  summand is non-positive with its null-space spanned by
  $\vect{\theta}\kronecker \vect{1}_N$, $\vect{\theta}\in\real^d$. On
  this null-space, the first summand becomes
    $-\alpha\vect{\theta}^\top\sum\nolimits_{i=1}^N(\nabla f^i(\vect{\theta} +
    \vect{\mathsf{x}}^\star)-\nabla f^i(\vect{\mathsf{x}}^\star))$, 
    which can only be zero when $\vect{\theta}=\vect{0}$, because
    $\sum_{i=1}^N \nabla f^i(\vect{\sf{x}})=\nabla f(\vect{\sf{x}})$
    and the global cost function is strictly convex by
    assumption. Then, the two summands of $ \dot{V}$ can be zero
    simultaneously only when $\vect{y}^i=\vect{0}$, for all
    $i\until{N}$, which is equivalent to $\vect{z}=\vect{0}$. Thus,
    $\dot{V} $ is negative semi-definite, with
    $\dot{V}(\vect{z},\vect{w}_{2:N})=0$ happening on the set
    $\mathcal{S} = \{(\vect{z},\vect{w}_{2:N}) \in \real^{Nd} \times
    \real^{(N-1)d} \;|\; \vect{z}=\vect{0}\}$.  Note
    that~\eqref{eq::Alg_Separated} on $\mathcal{S}$ reduces to
    $\dot{\vect{w}}_{2:N}=\vect{0}$, $\dvect{z}_1=\vect{0}$, and
    $\dvect{z}_{2:N}= -\vect{w}_{2:N}$.  Therefore, the only
    trajectory of~\eqref{eq::Alg_Separated} that remains in
    $\mathcal{S}$ is the equilibrium point
    $(\vect{z}_1=\vect{0},\vect{z}_{2:N}=\vect{0},\vect{w}_{2:N}=\vect{0})$.
    The LaSalle invariance principle (cf.~\cite[Theorem 4.4 and
    Corollary 4.2]{HKK:02}) now implies that the equilibrium is
    globally asymptotically stable or, in other words,
    $\vect{x}^i\to\vect{\mathsf{x}}^\star,~i\until{N}$ globally
    asymptotically.~\boxend
\end{pf}

\begin{rem}\longthmtitle{Simplification of~\eqref{eq::Alg}}
  {\rm 
    For strictly convex local cost functions, using the LaSalle
    function of the proof of Theorem~\ref{thm::UD_Conv_Alg2}, one can
    show that the asymptotic converge to the optimizer, starting from
    the initial condition stated in~Theorem~\ref{thm::UD_Conv_Alg2},
    is also guaranteed for the following algorithm over connected
    graphs \vspace{-0.1in}
     \begin{subequations}
       \begin{align*}
         \dot{\vect{v}}^i & = \sum\nolimits_{j=1}^N
         \mathsf{a}_{ij}(\vect{x}^i-\vect{x}^j),
         \\
         \dot{\vect{x}}^i & =-\nabla f^i(\vect{x}^i)-\vect{v}^i . \eqoprocend
     \end{align*}
   \end{subequations}
 }
 \end{rem}

\section{Continuous-time Evolution with Discrete-Time
  Communication}\label{sec::Discrete}

Here, we investigate the design of continuous-time coordination
algorithms with discrete-time communication to solve the distributed
optimization problem of Section~\ref{se:problem-statement}.  The
implementation of~\eqref{eq::Alg} requires continuous-time
communication among the agents. While this abstraction is useful for
analysis, in practical scenarios communication is only available at
discrete instants of time. This observation motivates our study here.
Throughout the section, we deal with communication topologies
described by connected graphs. Our results build on the discussion of
Section~\ref{sec::Alg}, particularly the identification of Lyapunov
functions for stability.

We start by introducing some useful conventions.  At any given time $t
\in \realnonnegative$, let $\Hvect{x}^j$ be the last known state of
agent $j \until{N}$ transmitted to its in-neighbors.  If
$\{t^i_{k}\} \subset \realnonnegative$ denotes the times at which
agent $i$ communicates with its in-neighbors, then one has
$\Hvect{x}^i = \vect{x}^i(t^i_{k})$ for
$t\in[t^i_{k},t^i_{k+1})$.  Consider the next implementation of
the algorithm~\eqref{eq::Alg} with discrete-time communication,
\vspace{-0.1in}
\begin{subequations}\label{eq::Alg-EvnTrig}
  \begin{align}
    \dvect{v}^i & =\alpha\beta\sum\nolimits_{j=1}^N
    \mathsf{a}_{ij}(\Hvect{x}^i-\Hvect{x}^j), \label{eq::Alg-EvnTrig-a}
    \\
    \dvect{x}^i & =-\alpha \nabla f^i(\vect{x}^i)-\beta\sum\nolimits_{j=1}^N
    \mathsf{a}_{ij}(\Hvect{x}^i-\Hvect{x}^j)-\vect{v}^i. 
    \label{eq::Alg-EvnTrig-b}
  \end{align}
\end{subequations}
Clearly, the evolution of~\eqref{eq::Alg-EvnTrig} depends on the
sequences of communication times for each agent. Here, we consider
three scenarios. Section~\ref{sec:periodic} studies periodic
communication schemes where all agents communicate synchronously at
$\Delta$ intervals of time, i.e., $t^i_{k}=t_k = \Delta k$ for all
$i\until{N}$. We provide a characterization of the periods that
guarantee the asymptotic convergence of~\eqref{eq::Alg-EvnTrig} to the
optimizer.  In general, periodic schemes might result in a wasteful
use of the communication resources because of the need to account for
worst-case situations in determining appropriate periods.  This
motivates our study in Section~\ref{sec:event-triggered} of
event-triggered communication schemes that tie the communication times
to the network state for greater efficiency. We discuss two
event-triggered communication implementations, a centralized
synchronous one and a distributed asynchronous one. In both cases, we
pay special attention to ruling out the presence of Zeno behavior (the
existence of an infinite number of updates in a finite interval of
time).

\subsection{Periodic Communication}\label{sec:periodic}

The following result provides an upper bound on the size of admissible
stepsizes for the execution of~\eqref{eq::Alg-EvnTrig} over connected
graphs with periodic communication schemes.

\begin{thm}[Convergence of~\eqref{eq::Alg-EvnTrig} with periodic
    communication]\label{thm::periodic}
    Let $\GG$ be a connected graph. Assume the local cost function
    $f^i$, $i\!\until{N}$, is $m^i$-strongly convex, differentiable,
    and its gradient is $M^i$-Lipschitz on~$\real^d$.  Given
    $\alpha,\beta>0$, consider an implementation
    of~\eqref{eq::Alg-EvnTrig} with agents communicating over $\GG$
    synchronously every $\Delta$ seconds starting at $t_0=0$, i.e.,
    $t^i_{k}=t_k = \Delta k$ for all $i\until{N}$.  Let $0<\eps<1$ and
    $\delta>0$ such that \vspace{-0.1in}
  \begin{align}\label{eq::trigger_error-phi}
    \phi = \frac{1}{2\mT} \MT^2 +
    \frac{1}{2\mT\alpha^2}\delta-1> 0,
  \end{align}
  where $\MT$ and $\mT$ are given in
  Theorem~\ref{thm::OpSlvrDirected}, and define \vspace{-0.1in}
  \begin{align}\label{eq::period_tau}
    \!  \tau\!= \!\frac{1}{\alpha
      \MT+1}\!\ln\!\Big(1+\frac{(\alpha
      \MT+1)\zeta}{\alpha \MT+1
      +\beta\lambda_N\sqrt{1+\alpha^2}(1+\zeta)}\Big),
  \end{align}
  where $\zeta^2 = \frac{2\eps{(1-\eps)\lambda_2 \min\{\delta,1\}}}{\alpha
    \beta\lambda_N^2\phi+ {4}\alpha^2\lambda_2(1+\phi)^2} $. Then, if
  $\Delta\in(0,\tau)$, the algorithm evolution starting from initial
  conditions $\vect{x}^i(0),\vect{v}^i(0)\in\reals^d$ with
  $\sum_{i=1}^N\vect{v}^i(0)=\vect{0}_d$ makes
  $\vect{x}^i(t)\to\vect{\mathsf{x}}^\star$ exponentially fast as
  $t\to\infty$, for all $i\until{N}$ with a rate of no less than $\frac{1}{4}\eps(\min\{\frac{1}{2},\delta\})/\bar{\lambda}_{\vect{E}}>0$.
\end{thm}
\begin{pf}
  We start by transferring the equilibrium point to the origin
  using~\eqref{eq::shift_orig} and then apply the change of
  variables~\eqref{eq::trans_stable} to write~\eqref{eq::Alg-EvnTrig}
  as
  \vspace{-0.1in}
  \begin{subequations}\label{eq::Alg-EvnTrig_Separated}
    \begin{align}
      \dvect{w}_{1} &= \vect{0}_d,
      \\
      \dvect{w}_{2:N} &
      =\alpha\beta(\rR^\top\lL\rR \kronecker
      \vect{\mathsf{I}}_d)(\vect{z}_{2:N}+\Tvect{z}_{2:N})
      \label{eq::Alg-EvnTrig_Separated-a},
      \\
      \dvect{z}_{1} &=- \alpha(\vect{\mathsf{r}}^\top\kronecker
      \vect{\mathsf{I}}_d)\vect{h},\label{eq::Alg-EvnTrig_Separated-b}
      \\
      \dvect{z}_{2:N}&=-\alpha(\rR^\top\kronecker
      \vect{\mathsf{I}}_d)\vect{h} \label{eq::Alg-EvnTrig_Separated-c}
      \\
      & \quad -\beta(\rR^\top\lL\rR \kronecker
      \vect{\mathsf{I}}_d)(\vect{z}_{2:N} +\Tvect{z}_{2:N})-\vect{w}_{2:N},
      \nonumber
    \end{align}
  \end{subequations}
  where $\Tvect{z}_{2:N}
  (t)\!=\!\vect{z}_{2:N}(t_k)\!-\!\vect{z}_{2:N} (t)$, for $t \in
  [t_k,t_{k+1})$, and $\vect{h}$ is given by~\eqref{eq::h}.  To study
  the stability
  of~\eqref{eq::Alg-EvnTrig_Separated-a}-\eqref{eq::Alg-EvnTrig_Separated-c},
  consider the Lyapunov function~\eqref{eq::LyapunovStronConvexUndi}
  with $\phi$ satisfying~\eqref{eq::trigger_error-phi}.  Its Lie
  derivative can be bounded by (details similar to the proof of
  Theorem~\ref{thm::UD_ConvRate_Alg} are omitted for~brevity)
  \vspace{-0.15in}
  \begin{align}\label{eq::V_dot_UAlg3}
    &\dot{V}  \leq -\frac{\delta{(1-\eps)} }{2}( \vect{z}_1^\top\vect{z}_1\! +\!
    \vect{z}_{2:N}^\top\vect{z}_{2:N}){-\frac{\delta\eps}{2}( \vect{z}_1^\top\vect{z}_1 \!+\!
    \vect{z}_{2:N}^\top\vect{z}_{2:N})}\nonumber
    \\
    &-\!\frac{\!(1\!-\!\eps)\!}{2} \vect{w}_{2:N}^\top\!\vect{w}_{2:N}\!-\!{\frac{\eps}{4} \vect{w}_{2:N}^\top\!\vect{w}_{2:N}}\!-\!\phi \alpha \beta \lambda_2(1\!-\!\eps)
    \vect{z}_{2:N}^\top\vect{z}_{2:N}\nonumber\\
    & \quad + \frac{1}{4\eps\lambda_2}
    \phi\alpha\beta\lambda_N^2{\Tvect{z}_{2:N}}^\top\Tvect{z}_{2:N} +
    \frac{1}{\eps}\alpha^2(\phi+1)^2{\Tvect{z}_{2:N}}^\top\Tvect{z}_{2:N}
    \nonumber
    \\
    & \quad\leq -\phi \alpha \beta \lambda_2(1-\eps)
    \vect{z}_{2:N}^\top\vect{z}_{2:N}-{\frac{1}{2}\eps\min\{\delta,\frac{1}{2}\}\vect{p}^\top\vect{p}}\nonumber
    \\
    & ~\quad -
    \frac{1}{2}{(1-\eps)}\min\{\delta,1\}\big(\vect{p}^\top\vect{p}-
    \zeta^{-2}
       {\Tvect{z}_{2:N}}^\top\Tvect{z}_{2:N}\big), \vspace{-0.10in}
  \end{align}
  where $\vect{p}=(\vect{z},\vect{w}_{2:N})$ and $\eps$ and $\zeta$
  are given in the theorem's statement. Observe that at each
  communication time $t_k$, $\|\Tvect{z}_{2:N} (t_k)\|=0$, then, it
  grows until next communication at time $t_{k+1}$ when it becomes
  zero again. Our proof proceeds by showing that if $t_{k+1}<
  t_{k}+\tau$, where~$\tau$ is given in~\eqref{eq::period_tau}, then
  we have the guarantee that
  \vspace{-0.1in}
  \begin{align}\label{eq::z-p}
    \|\Tvect{z}_{2:N} (t) \| < \zeta \, \|\vect{p} (t) \|, \quad
    t\in[t_k,t_{k+1}),
  \end{align}
  (note that, from~\eqref{eq::V_dot_UAlg3}, this guarantee ensures
  that $\dot{V}$ is negative definite for all $t\geq 0$). To this end,
  we study the dynamics of $q=\|\Tvect{z}_{2:N}\|/\|\vect{p}\|$ and
  find a lower bound on the time that it takes for $q$ to evolve from
  zero (recall $\Tvect{z}_{2:N}(t_k)=\vect{0}$) to~$\zeta$. Notice
  that \vspace{-0.1in}
  \begin{align*}
    & \dot{q} = \frac{\Tvect{z}_{2:N}^\top
      \overset{\textbf{.}}{\Tvect{z}}_{2:N}}{\|\Tvect{z}_{2:N}\|\|\vect{p}\|}
    - \frac{\|\Tvect{z}_{2:N}\|\vect{p}^\top\dvect{p}}{\|\vect{p}\|^3}
    \leq (1+q)\frac{\|\dvect{p}\|}{\|\vect{p}\|}\leq (1+q)\times
    \\
    &\frac{(\alpha \MT\!+\!\beta\lambda_N
      \sqrt{1\!+\!\alpha^2})\|\vect{p}\|\!+\!\|\vect{w}_{2:N}\|\!  +
      \!\beta\lambda_N\sqrt{1\!+\!\alpha^2}\|\Tvect{z}_{2:N}\|}{\|\vect{p}\|}
    \\
    & \leq(\alpha \MT+1)(1+q)+\beta\lambda_N\sqrt{1+\alpha^2}(1+q)^2.
  \end{align*}
  Here, we used in the first inequality $ d/dt (\Tvect{z}_{2:N})
  = - \dot {\vect{z}}_{2:N}$ and $\|\dvect{z}_{2:N}\|\leq
  \|\dvect{p}\|$ and in the second one the evolution of $\dvect{p}$
  given
  in~\eqref{eq::Alg-EvnTrig_Separated-a}-\eqref{eq::Alg-EvnTrig_Separated-c},
  $\|\rR^\top\lL\rR \|\leq \|\lL\| = \lambda_N$ and \vspace{-0.1in}
  \begin{align*}
    &\left\|
      \begin{bmatrix}
        -\beta\rR^\top\lL\rR 
        \\
        \alpha\beta\rR^\top\lL\rR  
      \end{bmatrix}
    \right\| = \beta\lambda_N\sqrt{1+\alpha^2}.
  \end{align*}
  Using the Comparison Lemma (cf.~\cite[Lemma 3.4]{HKK:02}), we
  conclude that $q(t,q_0)\leq \psi(t,\psi_0)$, where $\psi(t,\psi_0)$
  is the solution of $\dot{\psi}= (\alpha \MT+1)(1+\psi) +
  \beta\lambda_N\sqrt{1+\alpha^2}(1+\psi)^2$ satisfying
  $\psi(0,\psi_0)=\psi_0$. Then,
  \vspace{-0.1in}
  \begin{align*}
    & q(t,0)\leq \psi(t,0)
    \\
    & = \frac{ (\alpha
      \MT+1+\beta\lambda_N\sqrt{1+\alpha^2})(\text{e}^{\alpha
        \MT
        t+t}-1)}{-\beta\lambda_N\sqrt{1+\alpha^2}\text{e}^{\alpha
        \MT t+t}+\alpha \MT
      +1+\beta\lambda_N\sqrt{1+\alpha^2}}.
  \end{align*}
  The time $\tau$ for $\psi(\tau,0)\! =\!\zeta$ is given
  by~\eqref{eq::period_tau}.  Then, for $\{t_{k+1}-t_k\}_{k\in
    \integernonnegative}\!<\! \tau$, we have~\eqref{eq::z-p}, and as a
  result from~\eqref{eq::V_dot_UAlg3} we have
  $\dot{V}<\!-{\frac{1}{2}\eps\min\{\delta,\frac{1}{2}\}\vect{p}^\top\vect{p}}$.
  Thus, $\vect{z}\to \vect{0}$, as $t\to\infty$, which is equivalent
  to $\vect{x}^i\to\vect{\mathsf{x}}^\star$ as $t\to\infty$,
  exponentially fast with the rate given in the statement.
  \boxend
\end{pf}

\begin{rem}\longthmtitle{Dependence of the communication period on the
    design
    parameters}\label{rem:tau-depedence}
  {\rm The value of $\tau$ in Theorem~\ref{thm::periodic} depends on
    the graph topology, the parameters of the local cost functions,
    the design parameters $\alpha$ and $\beta$, and the variables
    $\eps$ and $\delta$.  One can use this dependency to maximize the
    value of~$\tau$.  Note that the argument of $\ln(.)$
    in~\eqref{eq::period_tau} is a monotonically increasing function
    of $\zeta>0$.  Therefore, the smaller the value of $\beta$, the
    larger the value of $\tau$. However, the dependency of $\tau$ on
    the rest of the parameters listed above is more complex.  For
    given local cost functions, fixed network topology and fixed
    values of $\alpha$, $\beta$, the maximum value of $\zeta$ is when
    $\phi+1$ is at its minimum and $\eps\lambda_2
    \min\{1-\eps,\delta\}$ is at its maximum.
    \oprocend }
\end{rem}

\subsection{Event-Triggered Communication}\label{sec:event-triggered}

This section studies the design of event-triggered communication
schemes for the execution of~\eqref{eq::Alg-EvnTrig}. In contrast to
periodic schemes, event-triggered implementations tie the
determination of the communication times to the current network state,
resulting in a more efficient use of the resources.  The proof of
Theorem~\ref{thm::periodic} reveals that the satisfaction of
condition~\eqref{eq::z-p} guarantees the monotonic evolution of the
Lyapunov function, which in turn ensures the correct asymptotic
behavior of the algorithm.  One could therefore specify when
communication should occur by determining the times when this
condition is not satisfied. There is, however, a serious drawback to
this approach: the evaluation of the condition~\eqref{eq::z-p}
requires the knowledge of the global minimizer
$\vect{\mathsf{x}}^\star$, which is of course not available. To see
this, note that \vspace{-0.1in}
\begin{subequations}\label{eq::z_2N-w_2:N_x_v}
  \begin{align}
    \|\Tvect{z}_{2:N}\| & =\|\PPi(\vect{x}(t_k)-\vect{x})\|
    , \label{eq::z_2N-w_2:N_x_v-a}
    \\
    \|\vect{p}\| & =
    \sqrt{\|\vect{x}-\Bvect{x}\|^2+\|\PPi(\vect{v}-\Bvect{v})\|^2}, 
    \label{eq::z_2N-w_2:N_x_v-b}
  \end{align}
\end{subequations}
where we have used 
\eqref{eq::Rrs} and~\eqref{eq::sum_v} (recall $\PPi\!=\!
\pPi_N\kronecker\vect{\mathsf{I}}_d$). From~\eqref{eq::v_bar},
$\Bvect{v}^i\!=\!-\alpha\nabla f^i(\vect{\mathsf{x}}^\star)$ for
$i\until{N}$, and thus the evaluation of the triggering
condition~\eqref{eq::z-p} requires knowledge of the global
optimizer. Our forthcoming discussion shows how one can circumvent
this problem. We first consider the design of centralized triggers
requiring global network knowledge and then discuss triggering schemes
that only rely on inter-neighbor~interaction.

\subsubsection{Centralized Synchronous Implementation}

Here, we present a centralized event-triggered scheme to
determine the sequence of synchronous communication times
in~\eqref{eq::Alg-EvnTrig}. Our discussion builds upon the examination
of the Lie derivative of the Lyapunov function used in the proof of
Theorem~\ref{thm::periodic} and the observations made above regarding
the lack of knowledge of the global optimizer.
From~\eqref{eq::z_2N-w_2:N_x_v}, we see that an event-triggered law
should not employ $\vect{p}$, but rather rely on $\Tvect{z}_{2:N}$
and~$\vect{z}_{2:N}$, to be independent of
$\vect{\mathsf{x}}^\star$. With this in mind, the examination of the
upper bound~\eqref{eq::V_dot_UAlg3} on~$\dot{V}$ reveals that, if
\vspace{-0.1in}
\begin{align}\label{eq::triggerLawUndirected}
  \|\Tvect{z}_{2:N} (t)\|^2 & = \|\PPi(\vect{x}(t_k)-\vect{x} (t))\|^2
  \notag
  \\
  & \leq \kappa \|\vect{z}_{2:N} (t) \|^2 =
  \kappa\|\PPi\vect{x}(t)\|^2,
\end{align}
where $\kappa$ is shorthand notation for
\begin{align}\label{eq::kappa}
  \kappa = 2 \frac{\eps \delta\lambda_2 + 2 \phi\alpha\beta\lambda_2^2
    \eps^2(1-\eps)}{\alpha\beta\phi\lambda_N^2 +2\lambda_2
    \alpha^2(1+\phi)^2},
\end{align}
(here $0<\eps<1$ and $\phi$ is given
by~\eqref{eq::trigger_error-phi}), then we have
\begin{align}\label{eq::CentEvenResult_Vdot}
  \dot{V} \!&\! \leq \!\! - \phi\alpha\beta\lambda_2(1\!-\!\eps)^2
  \vect{z}_{2:N}\!^\top\!  \vect{z}_{2:N}\!-\! \frac{1}{2}(1\!-\!\eps)
  \vect{w}_{2:N}\!^\top\!\vect{w}_{2:N} \!\!
  \\
  & \quad \!\!\!-\frac{1}{2}\delta \vect{z}_1\!^\top\!\vect{z}_1\!
  \leq\!-\frac{1}{2}\min\{\delta,\!2\phi\alpha\beta\lambda_2(1-\eps)^2,\!(1\!-\!\eps)\}\vect{p}\!^\top\!\vect{p}
  . \nonumber
\end{align}
Then, we can reproduce the proof of Theorem~\ref{thm::periodic} and
conclude the exponential convergence to the optimal solution.
Accordingly, the sequence of synchronous communication times
$\{t_{k}\}_{k \in \integernonnegative} \subset \realnonnegative$
for~\eqref{eq::Alg-EvnTrig} should be determined
by~\eqref{eq::triggerLawUndirected}.  However, for a truly
implementable law, one should rule out Zeno behavior, i.e., the
sequence of times does not have any finite accumulation point.
However, observing~\eqref{eq::triggerLawUndirected}, one can see that
Zeno behavior will arise at least near the agreement surface $\PPi
\vect{x}=\vect{0}_{dN}$.  The next result details how we address this
problem. 

\begin{thm}[Convergence of~\eqref{eq::Alg-EvnTrig} with Zeno-free
  centralized event-triggered
  communication]\label{thm::UD_EventTriggered_Alg3}
  Let $\GG$ be a connected graph.  Assume the local cost function
  $f^i$, $i\!\until{N}$, is $m^i$-strongly convex, differentiable, and
  its gradient is $M^i$-Lipschitz on~$\real^d$.  Consider an
  implementation of~\eqref{eq::Alg-EvnTrig} with agents communicating
  over $\GG$ synchronously at $\{t_{k}\}_{k \in \integernonnegative}
  \subset \realnonnegative$, starting at $t_0=0$,
  \vspace{-0.1in}
  \begin{multline}\label{eq::triggerLawUndirected_aug}
    t_{k+1} = \argmax \setdef{t \in [t_k+\tau,\infty)}
    {
          \\
      \|\PPi (\vect{x}(t_k)-\vect{x}(t))\|^2 \le \kappa\| \PPi
      \vect{x}(t)\|^2 },
  \end{multline}
  where $\tau$ and $\kappa<1$ are defined in~\eqref{eq::period_tau}
  and~\eqref{eq::kappa}, respectively.  Then, for any given
  $\alpha,\beta>0$ and each $i\until{N}$, the algorithm evolution
  starting from initial conditions
  $\vect{x}^i(0),\vect{v}^i(0)\in\reals^d$ with
  $\sum_{i=1}^N\vect{v}^i(0) = \vect{0}_d$ makes $\vect{x}^i(t) \to
  \vect{\mathsf{x}}^\star$ exponentially fast as $t\to\infty$ with a
  rate no less than $\frac{1}{4}(\min\{\delta,\!2\phi\alpha\beta
  \lambda_2(1-\eps)^2,\!(1\!-\!\eps),
  \frac{1}{2}\eps\})/\bar{\lambda}_{\vect{E}}>0$.
\end{thm}
\begin{pf}
  We first show
  $\kappa<1$.  This is an important property guaranteeing that, if
  agents start in agreement at a point other than the
  optimizer~$\vect{\mathsf{x}}^\star$, then the
  condition~\eqref{eq::triggerLawUndirected} is eventually violated,
  enforcing information updates.  Notice that (a) $ 4
  \eps^2(1-\eps)\lambda_2^2<\lambda_N^2$ and (b) $\eps
  \delta<\alpha^2(1+\phi)^2$ imply that the numerator
  in~\eqref{eq::kappa} is smaller than its denominator, and hence
  $\kappa<1$. (a) follows from noting that the maximum of $4
  \eps^2(1-\eps)$ for $\eps \in (0,1)$ is $16/27<1$ and the fact that
  $\lambda_2 \le \lambda_N$. We prove (b) reasoning by contradiction.
  Assume $ \delta>\alpha^2(1+\phi)^2$ or equivalently
  $\alpha(1+\phi)-\sqrt{\delta}<0$.
  Using~\eqref{eq::trigger_error-phi} and multiplying both sides of
  the inequality by $2\mT\alpha$, we obtain
    \begin{equation*}
    \alpha^2\MT^2+\delta-2\alpha \mT\sqrt{\delta}
   \! =\! (\sqrt{\delta} - \alpha \mT)^2 \!+\!
    \alpha^2(\MT^2-\mT^2)\!<\!0,
  \end{equation*}
  which, since $\MT\geq \mT$, is a contradiction. Having established
  the consistency of~\eqref{eq::triggerLawUndirected}, consider now
  the candidate Lyapunov function $V$
  in~\eqref{eq::LyapunovStronConvexUndi} and let $t_k$ be the last
  time at which a communication among all neighboring agents occurred.
  From the proof of Theorem~\ref{thm::periodic}, we know that the time
  derivative of $V$ is negative,
  $\dot{V}<-{\frac{1}{2}\eps\min\{\delta,\frac{1}{2}\}\vect{p}^\top\vect{p}}$
  as long as $t < t_k+\tau$.  After this
  time,~\eqref{eq::CentEvenResult_Vdot} shows that as long
  as~\eqref{eq::triggerLawUndirected} is satisfied, $\dot{V}$ is
  negative, and exponential convergence follows. \boxend
\end{pf}

Interestingly, given that~\eqref{eq::triggerLawUndirected} does not
use the full state of the network but instead relies on the
disagreement, one can interpret it as an output feedback
event-triggered controller.  Guaranteeing the existence of lower
bounded inter-execution times for such controllers is in general a
difficult problem, see e.g.,~\citep{MCFD-WPMHH:12}.
Augmenting~\eqref{eq::triggerLawUndirected} with the condition
$t_{k+1} \ge t_k + \tau$ results in Zeno-free executions by lower
bounding the inter-event times by~$\tau$. The knowledge of this value
also allows the designer to compute bounds on the maximum energy spent
by the network on communication.

\subsubsection{Distributed Asynchronous Implementation}

We present a distributed event-triggered scheme for determining the
sequence of communication times in~\eqref{eq::Alg-EvnTrig}.  At each
agent, the execution of the communication scheme depends only on local
variables and the triggered states received from its neighbors. This
naturally results in asynchronous communication. We also show that the
resulting executions are free from Zeno behavior.

\begin{thm}[Convergence of~\eqref{eq::Alg-EvnTrig} with Zeno-free
  distributed event-triggered
  communication]\label{thm::UD_DecentralziedEventTriggered}
  Let $\GG$ be a connected graph. Assume $f^i$, $i\!\until{N}$, is
  $m^i$-strongly convex, differentiable, and its gradient is
  $M^i$-Lipschitz on~$\real^d$.  For $\vect{\eps} \in
  \realpositive^N$, consider an implementation
  of~\eqref{eq::Alg-EvnTrig} where agent $i \until{N}$ communicates
  with its neighbors in~$\GG$ at times $\{t^i_{k}\}_{k \in
    \bar{Z}^i\subseteq \integernonnegative} \subset \realnonnegative$,
  starting at $t^i_0=0$,\vspace{-0.1in}
  \begin{align}\label{eq::TrigLaw_Distributed}
    & t^i_{k+1} = \argmax \setdef{t \in [t^i_{k},\infty)} {
      \\
      & 4
      {\mathsf{d}}_{\text{out}}^i\|\Hvect{x}^i(t)\!-\!\vect{x}^i(t)\|^2\leq
      \sum\nolimits_{j=1}^N\!\!\!
      {\mathsf{a}}_{ij}\|\Hvect{x}^i(t)-\Hvect{x}^j(t)\|^2\! +\!
      (\eps^i)^2 }. \nonumber
  \end{align}
  Given $\alpha>0$, let $\beta,\phi>0$ satisfy $ \phi+1>
  4\MT$ and
  \begin{align}\label{eq::another-check}
   \!\!\!\gamma'\!=\! \alpha^2\!(\phi\!+\!1)\mT\! +\!
    \tfrac{9}{2}\beta\Hlambda_2\phi\alpha \!-\!  4\alpha^2(
    \MT\mT\!+\!(\phi\!+\!1)\!^2) \!\!>\!0.
  \end{align}
  Then, for each $i \until{N}$, the evolution starting from initial
  conditions $\vect{x}^i(0),\vect{v}^i(0)\in\reals^d$ with
  $\sum_{i=1}^N\vect{v}^i(0) = \vect{0}_d$ makes $\|\vect{x}^i(t) -
  \vect{\mathsf{x}}^\star\| \leq
  \frac{\phi\alpha\beta\bar{\lambda}_F}{4\eta\underline{\lambda}_F}\|
  \vect{\eps}\|^2$ as $t\to\infty$ exponentially fast with a rate no
  less than $\eta/\bar{\lambda}_F$.  (Here, $\eta = \min\{
  \frac{7}{16},\frac{1}{9}\gamma'\}$, and $\underline{\lambda}_F$ and
  $\bar{\lambda}_F$ are the minimum and maximum eigenvalues of
  $\vect{F}$ in~\eqref{eq:matrix-F}).  Moreover, the inter-execution
  times of~$i$ are lower bounded by \vspace{-0.1in}
  \begin{equation}\label{eq::dist_event_time_lw_bnd}
    \tau^i = \frac{1}{\alpha M^i} \ln \Big(1+\frac{\alpha
      M^i\eps^i}{2\sqrt{ \mathsf{d}_{\text{out}}^i} (\alpha M^i+2\beta
      \mathsf{d}_{\text{out}}^i+1) \theta} \Big) ,
  \end{equation}  
  where $\theta = \frac{\bar{\lambda}_F}{\underline{\lambda}_F}
  \sqrt{\|\vect{x}(0)-\Bvect{x}\|^2 \!+\!\|\vect{v}(0) -
    \Bvect{v}\|^2} + \frac{\phi\alpha
    \beta\bar{\lambda}_F}{4\eta\underline{\lambda}_F}
  \|\vect{\eps}\|^2$.
\end{thm}
\begin{pf}
  Given an initial condition, let $[0,T)$ be the maximal interval on
  which there is no accumulation point in $\{t_k\}_{k \in \bar{Z}} =
  \cup_{i=1}^N \cup_{k \in \bar{Z}^i\subseteq\integernonnegative}
  t^i_{k} $.  Note that $T>0$, since the number of agents is finite
  and, for each $i\until{N}$, $\eps^i>0$ and $\Tvect{x}^i(0) =
  \Hvect{x}^i(0)-\vect{x}^i(0)=\vect{0}$.  The
  dynamics~\eqref{eq::Alg-EvnTrig}, under the event-triggered
  communication scheme~\eqref{eq::TrigLaw_Distributed}, has a unique
  solution in the time interval $[0,T)$.  Next, we use Lyapunov
  analysis to show that the trajectory is bounded during~$[0,T)$.
  Consider the function $V$ given in~\eqref{eq::Lya_Digraph}, whose
  Lie derivative
  along~\eqref{eq::Alg-EvnTrig_Separated-a}-\eqref{eq::Alg-EvnTrig_Separated-c}
  is \vspace{-0.1in}
  \begin{align*}
    \dot{V} & = -\frac{1}{9}\alpha^2(\phi+1)\vect{y}^\top\vect{h} -
    \frac{7}{16}\vect{w}_{2:N}^\top\vect{w}_{2:N}
    \\
    & \quad
    -\|\frac{3}{4}\vect{w}_{2:N}+\frac{2\alpha}{3}(\rR^\top\kronecker
    \vect{\mathsf{I}}_d)\vect{h}+\frac{2\alpha}{3}(\phi+1)\vect{z}_{2:N})\|^2
    \\
    & \quad +\frac{4}{9}\alpha^2\|(\rR^\top\kronecker
    \vect{\mathsf{I}}_d)\vect{h}\|^2 +
    \frac{4}{9}\alpha^2(1+\phi)^2\vect{z}_{2:N}^\top\vect{z}_{2:N}
    \\
    & \quad -\frac{\phi\alpha\beta}{2} \vect{z}_{2:N}^\top \!(\rR^\top\lL\rR
    \!\kronecker\!  \vect{\mathsf{I}}_d) \vect{z}_{2:N}+\frac{\phi\alpha\beta}{2}s,
  \end{align*}
  where $s = - \vect{z}_{2:N}^\top(\rR^\top \lL\rR
    \!\kronecker\!  \vect{\mathsf{I}}_d)
    \vect{z}_{2:N}-2 \vect{z}_{2:N}^\top(\rR^\top \lL\rR
    \!\kronecker\! \vect{\mathsf{I}}_d)\Tvect{z}_{2:N}$, and $\Tvect{z}_{2:N}=\Hvect{z}_{2:N}-\vect{z}_{2:N}$. Using the assumptions on the local cost
  functions and following steps similar to those taken in the proof of
  Theorem~\ref{thm::OpSlvrDirected} to lower bound
  $\vect{y}^\top\vect{h}$ and upper bound $\|(\rR^\top\kronecker
  \vect{\mathsf{I}}_d)\vect{h}\|$ along with using $\lambda_2 \vect{z}_{2:N}^\top  \vect{z}_{2:N}\leq \vect{z}_{2:N}^\top \!(\rR^\top\lL\rR
  \!\kronecker\!  \vect{\mathsf{I}}_d) \vect{z}_{2:N}$ , one can show that
  \vspace{-0.1in}
  \begin{align*}
    \dot{V} &\leq -\eta\|\vect{p}\|^2+
    \frac{\phi\alpha \beta}{2}s ,
  \end{align*}
  where 
  $\vect{p}=(\vect{z},\vect{w}_{2:N})$.  
  Next, we show $s \leq \frac{1}{2}\|\vect{\eps}\|^2$ for $t \in
  [0,T)$.  Using $\rR \rR^\top=\pPi_N$, $\lL\pPi_N=\pPi_N\lL=\lL$,
  $\vect{z}_{2:N} = (\rR^\top\kronecker
  \vect{\mathsf{I}}_d)\vect{y}=(\rR^\top\kronecker
  \vect{\mathsf{I}}_d)\vect{x}$, and $\Tvect{z}_{2:N} =
  (\rR^\top\kronecker
  \vect{\mathsf{I}}_d)\Tvect{y}=(\rR^\top\kronecker
  \vect{\mathsf{I}}_d)\Tvect{x}$, we get
  $s=-\vect{x}^\top\LL\vect{x}-2\vect{x}^\top\LL\Tvect{x}$.   Then, \vspace{-0.1in}
  \begin{align*}
    s &=-(\Hvect{x}-\Tvect{x})^\top \LL
    (\Hvect{x}-\Tvect{x})-2(\Hvect{x}-\Tvect{x})^\top\LL\Tvect{x}
    \\
    &=\Tvect{x}^\top\LL\Tvect{x}-\Hvect{x}^\top\LL\Hvect{x}.
  \end{align*}
  Given $\LL=\vect{\mathsf{D}}_{\text{out}}-\vect{\mathsf{A}}$ and 
  $\vect{\mathsf{D}}_{\text{out}}+\vect{\mathsf{A}}\geq 0$, we have
  \vspace{-0.1in}
  \begin{align*}
    \Tvect{x}^\top\LL\Tvect{x} &
     \leq 2\Tvect{x}^\top(\vect{\mathsf{D}}_{\text{out}}\kronecker
    \vect{\mathsf{I}}_d)\Tvect{x} =
    2\sum\nolimits_{i=1}^N\mathsf{d}_{\text{out}}^i\|\Tvect{x}^i\|^2.
  \end{align*}
  Therefore, we can write (recall
  $\Tvect{x}^i=\Hvect{x}^i-\vect{x}^i$) \vspace{-0.1in}
  \begin{align*}
    s
    & = \frac{1}{2}\sum\nolimits_{i=1}^N
    \big(4\mathsf{d}_{\text{out}}^i\|\Hvect{x}^i-\vect{x}^i\|^2
    -\sum\nolimits_{j=1}^N \mathsf{a}_{ij}\|\Hvect{x}^i-\Hvect{x}^j\|^2\big),
  \end{align*}
  which, with~\eqref{eq::TrigLaw_Distributed}, yields $s \!\leq\!
  \frac{1}{2}\|\vect{\eps}\|^2$ for $t \!\in\! [0,T)$.~Then,
  \vspace{-0.1in}
  \begin{align*}
    \dot{V} 
        & \leq -\eta\|\vect{p}\|^2+\frac{\phi\alpha
      \beta}{4}\|\vect{\eps}\|^2,\quad t\in[0,T),
  \end{align*}
   Recall from the proof of
  Theorem~\ref{thm::OpSlvrDirected} that
  $\underline{\lambda}_F\|\vect{p}\|^2\leq V (\vect{p})\leq
  \bar{\lambda}_F\|\vect{p}\|^2$. Then, using the Comparison Lemma
  (cf.~\cite[Lemma 3.4]{HKK:02}), we deduce that
  \vspace{-0.1in}
  \begin{align}\label{eq::trajec_bound_decent_event}
    \!\!\|\vect{p}(t)\| & \leq\! \frac{1}{\underline{\lambda}_F}
    \|V(0)\|\e^{-\frac{\eta}{\bar{\lambda}_F}t}\!+
    \frac{\phi\alpha\beta\bar{\lambda}_F
      \|\vect{\eps}\|^2}{4\eta\underline{\lambda}_F}(1\!-\!
    \e^{-\frac{\eta}{\bar{\lambda}_F}t})
    \nonumber\\
    &\!\!\!\leq \!\frac{\bar{\lambda}_F}{\underline{\lambda}_F}
    \|\vect{p}(0)\|\!\e^{-\frac{\eta}{\bar{\lambda}_F}t}\!\!+
    \frac{\phi\alpha\beta\bar{\lambda}_F\|
      \vect{\eps}\|^2}{4\eta\underline{\lambda}_F}
    (1\!-\!\e^{-\frac{\eta}{\bar{\lambda}_F}t}) , 
  \end{align}
  for $t\in [0,T)$.  Notice that regardless of value of $T$,
  \vspace{-0.1in}
  \begin{align}\label{eq::bound_p_0T}
    &\|\vect{p}(t) \|\leq
    \frac{\bar{\lambda}_F}{\underline{\lambda}_F}\|\vect{p}(0)\| +
    \frac{\phi\alpha\beta\bar{\lambda}_F}{4\eta
      \underline{\lambda}_F}\|\vect{\eps}\|^2,
  \end{align} 
  for $t\in [0,T)$.  Notice that the right-hand side corresponds
  to~$\theta$. This can be seen by noting that,
  from~\eqref{eq::trans_stable}, we have $\|\vect{z}\| =
  \|\vect{x}-\Bvect{x}\|$ and $\|\vect{w}\| =
  \|\vect{w}_{2:N}\|=\|\vect{v}- \Bvect{v}\|$ (recall
  $\sum_{i=1}^N\vect{v}^i(0)=\vect{0}$ results in
  $\vect{w}_1=\vect{0}$ for all $t\geq0$).

  Our final objective is to show that $T=\infty$.  To achieve this, we
  firt establis a lower bound on the inter-execution times of any
  agent. To do this, we determine a lower bound on the time it takes
  $i \until{N}$ to have $\|\Hvect{x}^i-\vect{x}^i\|$ evolve from $0$
  to $\eps^i/(2\sqrt{ {\mathsf{d}}_{\text{out}}^i})$.
  Using~\eqref{eq::v_bar} and~\eqref{eq::Alg-EvnTrig-b}, we have
  \vspace{-0.1in}
  \begin{align*}
    &\frac{d}{dt} \|\Hvect{x}^i-\vect{x}^i\|=
    -\frac{(\Hvect{x}^i-\vect{x}^i)^\top\dvect{x}^i}{\|\Hvect{x}^i-\vect{x}^i\|}
    \leq \|\dvect{x}^i\|
    \\
    & = \|-\alpha(\nabla f^i(\vect{x}^i)-\nabla
    f^i(\vect{\mathsf{x}}^{\star})) - \beta\sum\nolimits_{j=1}^N
    \mathsf{a}_{ij}(\Hvect{x}^i\!-\!\Hvect{x}^j)
    \\
    & \quad -\!(\vect{v}^i+\alpha\nabla f^i(\vect{\mathsf{x}}^{\star})\|
    \\
    & \leq \alpha M^i \|\vect{x}^i-\vect{\mathsf{x}}^{\star}\| +
    \beta\sum\nolimits_{j=1}^N \mathsf{a}_{ij}\|\Hvect{x}^i\!-\!\Hvect{x}^j\|+
    \|\vect{v}^i- \Bvect{v}^i\|,
  \end{align*}
  which leads to
  \vspace{-0.1in}
  \begin{align*}
    & \frac{d}{dt} \|\Hvect{x}^i-\vect{x}^i\| \leq \alpha M^i
    \|\Hvect{x}^i-\vect{x}^i\|+ \alpha
    M^i\|\Hvect{x}^i-\vect{\mathsf{x}}^{\star}\| +
    \\
    & \quad \beta\sum\nolimits_{j=1}^N \mathsf{a}_{ij}(\|\Hvect{x}^i -
    \vect{\mathsf{x}}^{\star}\|+\|\Hvect{x}^j-\vect{\mathsf{x}}^{\star}\|)+
    \|\vect{v}^i- \Bvect{v}^i\|.
  \end{align*}
  From~\eqref{eq::bound_p_0T}, we have $\|\vect{x}(t)-\Bvect{x}\|\leq
  \theta$ and $\|\vect{v}(t)-\Bvect{v}\|\leq \theta$. This implies
  $\|\vect{v}^i(t)-\Bvect{v}^i\|\leq \theta$,
  $\|\Hvect{x}^i-\vect{\mathsf{x}}^{\star}\|\leq \theta$, and
  $\sum\nolimits_{j=1}^N \mathsf{a}_{ij}(\|\Hvect{x}^i -
  \vect{\mathsf{x}}^{\star}\| +
  \|\Hvect{x}^j-\vect{\mathsf{x}}^{\star}\|)\leq 2
  \mathsf{d}_{\text{out}}^i \theta$, for all $i\until{N}$.  Therefore,
  from the inequality above, $ \frac{d}{dt} \|\Hvect{x}^i-\vect{x}^i\|
  \leq \alpha M^i \|\Hvect{x}^i-\vect{x}^i\|+ c^i$, where $c^i =
  (\alpha M^i+2\beta \mathsf{d}_{\text{out}}^i+1) \theta$. Using the
  Comparison Lemma (cf.~\cite[Lemma 3.4]{HKK:02}) and the fact that
  $\|\Hvect{x}^i-\vect{x}^i(t^i_{k}) \|=0$, we deduce \vspace{-0.1in}
  \begin{align*}
    & \|\Hvect{x}^i-\vect{x}^i (t)\| \leq c^i 
(\e^{\alpha M^i (t-t^i_{k})}-1)/(\alpha
      M^i), \quad t \ge t^i_{k} . \vspace{-0.1in}
  \end{align*}
  Then, the time it takes $\|\Hvect{x}^i-\vect{x}^i\|$ to reach $
  \eps^i/(2\sqrt{ {\mathsf{d}}_{\text{out}}^i})$ is lower bounded by
  $\tau^i>0$ given by~\eqref{eq::dist_event_time_lw_bnd}.  To show
  $T=\infty$, we proceed by contradiction. 
  Suppose that $T<\infty$. Then, the sequence of events
  $\{t_k\}_{k\in\bar{Z}}$ has an accumulation point at $T$. Because we
  have a finite number of agents, this means that there must be an
  agent $i\until{N}$ for which $\{t^i_{k}\}_{k\in\bar{Z}^i}$ has an
  accumulation point at $T$, implying that agent $i$ transmits
  infinitely often in the time interval $[T-\Delta,T)$ for any
  $\Delta\in(0,T]$. However, this is in contradiction with the fact
  that inter-event times are lower bounded by $\tau^i>0$ on $[0,T)$.
  Having established $T=\infty$, note that this fact implies that
  under the event-triggered communication
  law~\eqref{eq::TrigLaw_Distributed}, the
  algorithm~\eqref{eq::Alg-EvnTrig} does not exhibit Zeno
  behavior. Furthermore, from \eqref{eq::trajec_bound_decent_event},
  we deduce that, for each $i \until{N}$, one has $\|\vect{x}^i(t)
  -\vect{\sf{x}}^{\star}\|\leq \|\vect{p} (t)\| \leq \frac{\phi\alpha
    \beta\bar{\lambda}_F}{4\eta\underline{\lambda}_F}\|\vect{\eps}\|^2$
  as $t\to\infty$, exponentially fast with a rate no worse than
  $\eta/\bar{\lambda}_F$.  \boxend
\end{pf}

Regarding the role of the design parameters and
condition~\eqref{eq::another-check}, we omit for space reasons
observations similar to the ones made in
Remark~\ref{rem:design-parameters}.  The lower bound on the
inter-event times allows the designer to compute bounds on the maximum
energy spent by each agent on communication during any given time
interval.  Since the total number of agents is finite and each agent's
inter-event times are lower bounded, it follows that the total number
of events in any finite time interval is finite. In general, an
explicit expression lower bounding the network inter-event times is
not available.  It is also worth noticing that the farther away the
agents start from the final convergence point (larger $\theta$
in~\eqref{eq::dist_event_time_lw_bnd}), the smaller the guaranteed
lower bound between inter-event times becomes.  As before, $\tau^i$
in~\eqref{eq::dist_event_time_lw_bnd} depends on the graph topology,
the parameters of the local cost function, the design parameters
$\alpha$ and $\beta$, and the variable $\eps^i$.  One can use this
dependency to maximize the value of~$\tau^i$ in a similar fashion as
discussed in Remark~\ref{rem:tau-depedence}.  For quadratic local cost
functions of the form $f^i(\vect{\mathsf{x}}) =
\frac{1}{2}(\vect{\mathsf{x}}^\top\vect{\mathsf{x}}+\vect{\mathsf{x}}^\top
\vect{a}^i+\vect{b}^i)$, $i\until{N}$ the claim of
Theorem~\ref{thm::UD_DecentralziedEventTriggered} holds for any
$\alpha,\beta>0$. Finally, we point out that the discrete-time
communication strategies introduced above enjoy similar privacy
preservation properties as the ones stated in
Proposition~\ref{prop::privLocalg}, but we omit the details here for
brevity.
\begin{figure*}[t!]
  \unitlength=0.5in \centering
  \psfrag*{x}[][cc][1.1]{\renewcommand{\arraystretch}{0.5}\begin{tabular}{c}$\ln(|x^i-\mathsf{x}^\star|)$\\~\end{tabular}}
\psfrag*{t}[][cc][1.1]{\renewcommand{\arraystretch}{2}\begin{tabular}{c}$t$\end{tabular}}
  \captionsetup[subfloat]{captionskip=-1pt}
  \subfloat[$\alpha=1$, $\beta=0.5$]{
    \includegraphics[width=2.2in, height=1.3in]{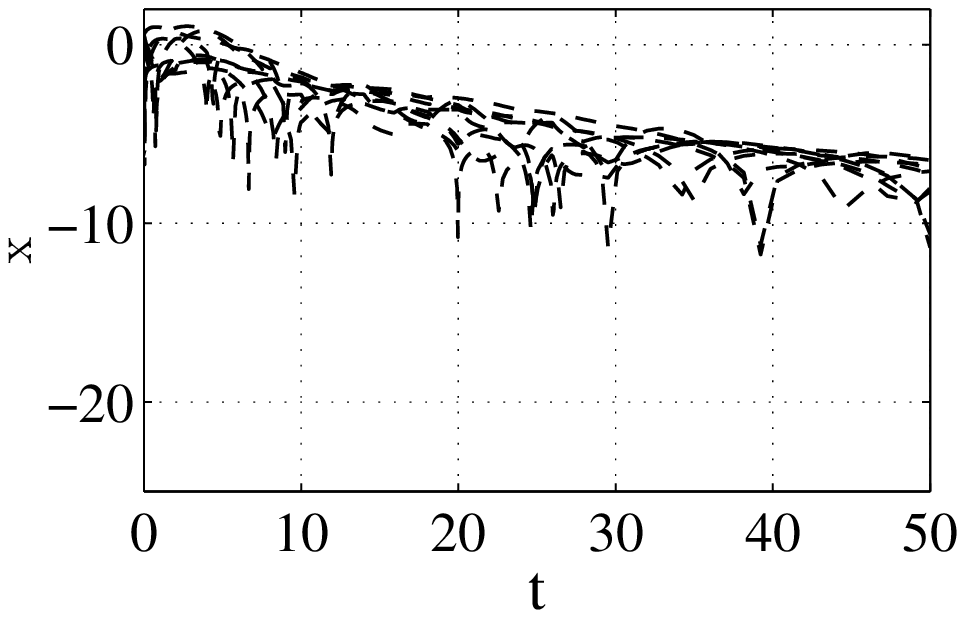}
  }
  \subfloat[$\alpha=1$, $\beta=1$]{
    \includegraphics[width=2.2in, height=1.3in]{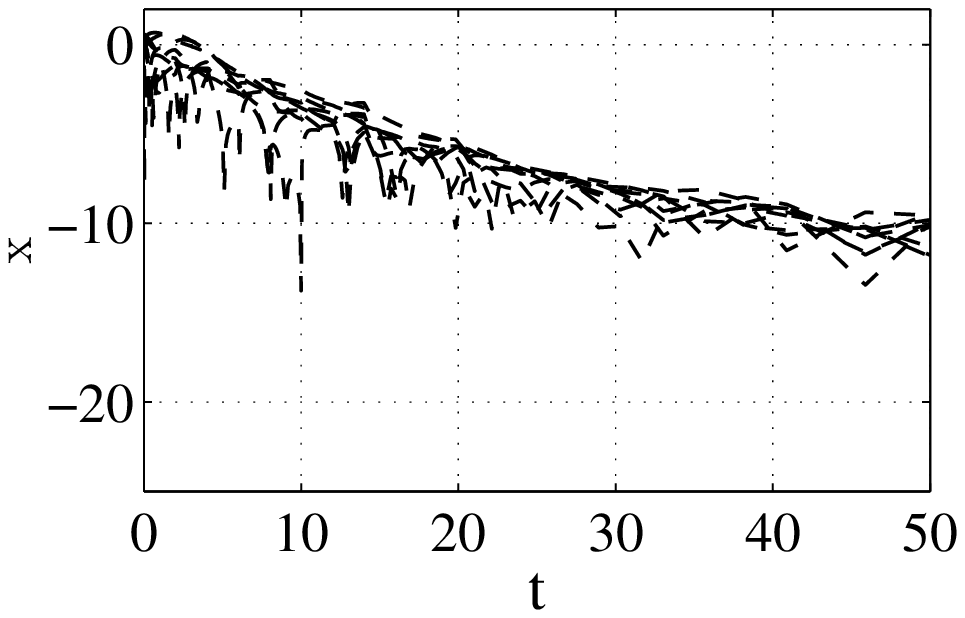} 
  }
  \subfloat[$\alpha=1$, $\beta=5$]{
    \includegraphics[width=2.2in, height=1.3in]{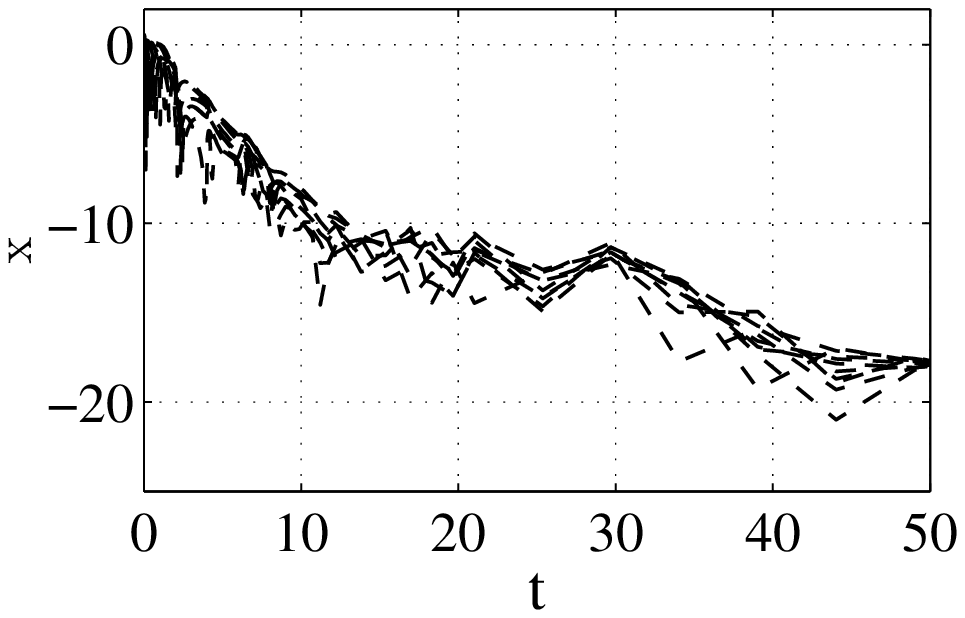} 
  }  \vspace{-2pt}
  \caption{Executions of~\eqref{eq::Alg} over a time-varying digraph
    that remains weight-balanced and strongly connected.
    }
\label{fig::Ex1sim}
\end{figure*} 

\section{Simulations}\label{sec::num}
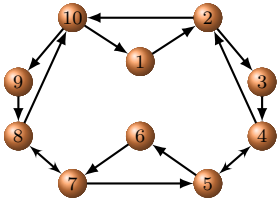
\begin{figure}[t]
  \centering
  \begin{tikzpicture}[auto,thick,scale=0.9, every
    node/.style={scale=0.9}]
    \tikzstyle{mynode}=%
    [%
    minimum size=12pt,%
    inner sep=0pt,%
    outer sep=0pt,%
    ball color=red!20!orange!70,
    shape=circle%
    ]
    \draw        
    (0,-1.9) node[mynode] (1) {{\scriptsize1}}
    (1,-1.25) node[mynode] (2) {{\scriptsize2}}
    (1.8,-2.2) node[mynode] (3) {{\scriptsize3}}
    (1.8,-3) node[mynode] (4) {{\scriptsize4}}
    (1,-3.7) node[mynode] (5) {{\scriptsize5}}
    (0,-3) node[mynode] (6) {{\scriptsize6}}
    (-1,-1.25) node[mynode] (10) {{\scriptsize10}}
    (-1.8,-2.2) node[mynode] (9) {{\scriptsize9}}
    (-1.8,-3) node[mynode] (8) {{\scriptsize8}}
    (-1,-3.7) node[mynode] (7) {{\scriptsize7}};
    \draw[-latex]  (1)->(2) ;
    \draw[-latex]  (2)->(3) ;
    \draw[-latex]  (3)->(4) ;
    \draw[latex'-latex'] (4) --  (5);
    \draw[-latex] (5)->(6) ;
    \draw[-latex]  (6)->(7) ;
    \draw[latex'-latex'] (7) --  (8);
    \draw[-latex] (9)->(8) ;
    \draw[-latex] (10)->(9) ;
    \draw[-latex] (10)->(1) ;
    \draw[-latex] (2)->(10) ;
    \draw[-latex] (4)->(2) ;
    \draw[-latex] (7)->(5) ;
    \draw[-latex] (8)->(10) ;
  \end{tikzpicture} \vspace{-2pt}
  \caption{One of the digraphs used in the simulations.
  }\label{fig::network}
\end{figure}

Here, we illustrate the performance of the algorithm~\eqref{eq::Alg}
and its implementation with discrete-time
communication~\eqref{eq::Alg-EvnTrig}. We consider a network of $10$
agents, with strongly convex local cost functions on $\real$ given by
\vspace{-0.1in}
\begin{alignat*}{2}
  f^1(x) & \! = \! 0.5\e^{-0.5 x} + 0.4 \e^{0.3x}, & \, f^2(x) & \! =
  \!(x-4)^2,
  \\
  f^3(x) & \! = \! 0.5x^2\ln(1+x^2)+x^2, & \, f^4(x) & \! =
  \!x^2+\e^{0.1x},
  \\
  f^5(x) & \! = \!\ln(\e^{-0.1x}\!+\!\e^{0.3x})+0.1x^2, & \, f^6(x) &
  \! = \!x^2\!/\!\ln(2\!+\!x^2),
  \\
  f^7(x) & \! = \! 0.2\e^{-0.2 x}+ 0.4 \e^{0.4x}, & \, f^8(x) &\! = \!
  x^4\!+\!2x^2\!+\!2 ,
  \\
  f^{9}(x) &\! = \!x^2/\sqrt{x^2+1}+0.1x^2 , & \, f^{10} (x) &\! =
  \!(x+2)^2.
\end{alignat*}
The gradient of the cost function of agents $1$, $4$, $7$, $8$ are
locally Lipschitz, while the rest are globally Lipschitz.
Figure~\ref{fig::Ex1sim} shows the executions of the
algorithm~\eqref{eq::Alg} for different values of $\beta$ when the
network topology alternates every $2$ seconds among three strongly
connected, weight-balanced digraphs (with unitary edge weights,
Figure~\ref{fig::network} shows one of these graphs).  Convergence is
achieved as guaranteed by Proposition~\ref{prop::OpSlvrDirectedSwtch}
(see also Remark~\ref{rem::Lips_compact}). The plot also shows that
larger values of $\beta$ result in faster convergence,
cf. Remark~\ref{rem:design-parameters}.
In all the simulations we ran, convergence is achieved for any
$\alpha,\beta>0$.
 
\begin{figure}[t]
  \unitlength=0.5in 
 \psfrag*{x}[][cc][1.1]{\renewcommand{\arraystretch}{0.5}\begin{tabular}{c}{\tiny$\ln(|x^i-\mathsf{x}^\star|)$}\\~\end{tabular}}
\psfrag*{t}[][cc][1.1]{\renewcommand{\arraystretch}{2}\begin{tabular}{c}{\tiny$t$}\end{tabular}}
 \captionsetup[subfloat]{captionskip=-1pt}
      \centering
   \subfloat[Algorithm~\eqref{eq::Alg-EvnTrig}: $\!\alpha\!=\!1,~~~~$ $\quad\quad\quad$ $~~~~~~~\!\beta\!=\!1$, $\Delta\!=\!0.5$~s.]{
    \includegraphics[height=.965in]{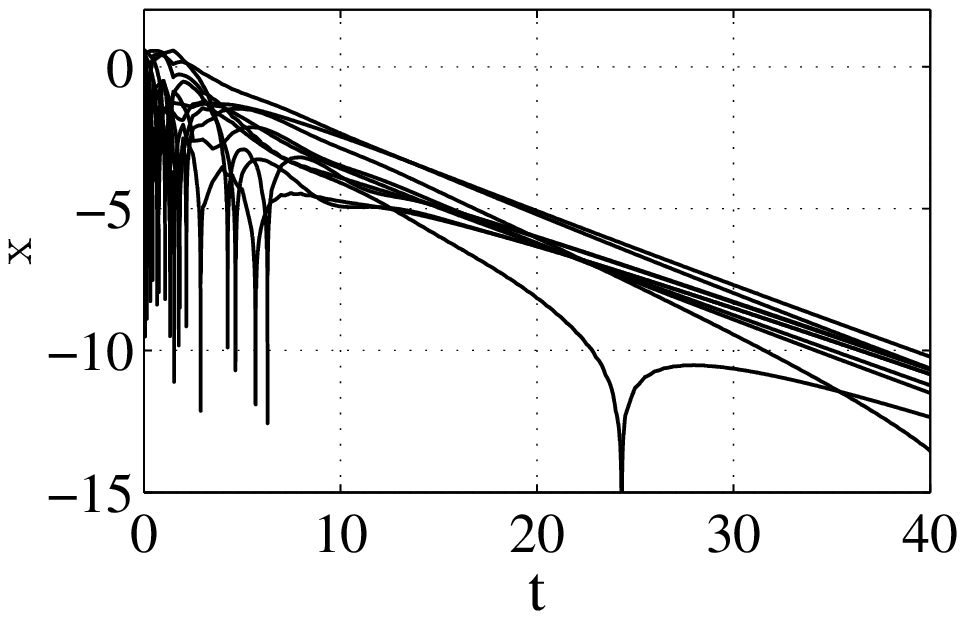}
  }\quad
  \subfloat[\!Algorithm~\eqref{eq::Alg-EvnTrig}: $\!\alpha\!=\!1,~~~~$ $\quad\quad\quad$ $~~~~~~\beta\!=\!0.5$, $\Delta\!=\!1$~s.]{
    \includegraphics[height=.965in]{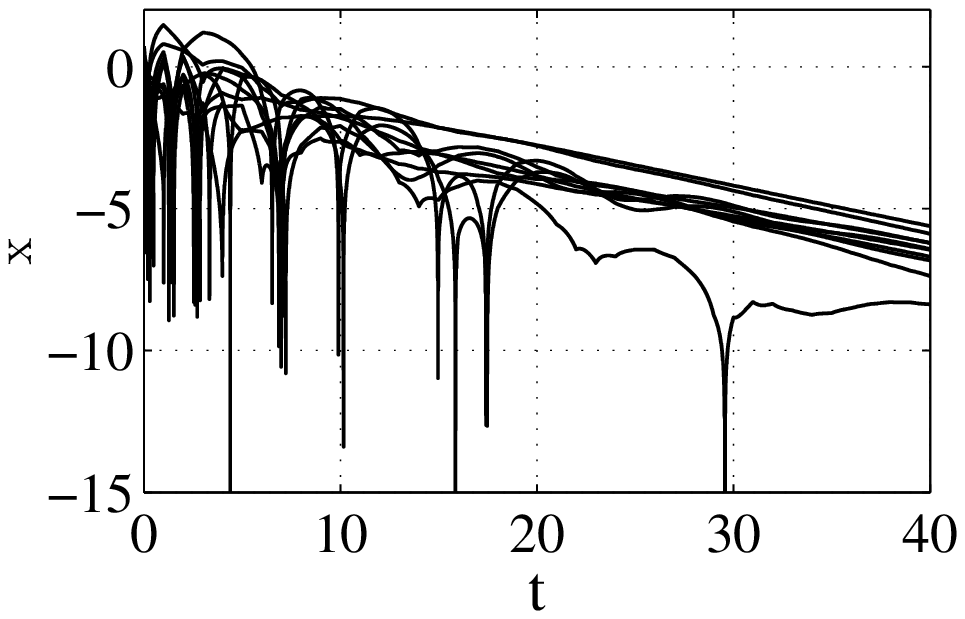} 
  } \vspace{-2pt}
  \caption{Performance evaluation of the
    algorithm~\eqref{eq::Alg-EvnTrig} when the communication is
    periodic.  
    }\label{fig::Ex2sim}
\end{figure}

Figures~\ref{fig::Ex2sim}(a)-(b) show executions
of~\eqref{eq::Alg-EvnTrig} with periodic communication over the
network depicted in~Figure~\ref{fig::network} for different $\beta$'s
and $\Delta$'s. Even though Theorem~\ref{thm::periodic} is established
for undirected graphs, our simulations show convergent behavior for
strongly connected and weight-balanced
digraphs. Figure~\ref{fig::Ex2sim} suggests a trade-off where larger
$\Delta$ (corresponding to smaller $\beta$, see
Remark~\ref{rem:tau-depedence}) result on savings on the energy
consumed by agents for communication at the cost of slower
convergence.

\begin{figure}[t!]
  \unitlength=0.5in \centering
  \psfrag*{x}[][cc][1.1]{\renewcommand{\arraystretch}{0.5}\begin{tabular}{c}{\tiny$\ln(|x^i-\mathsf{x}^\star|)$}\\~\end{tabular}}
  \psfrag*{t}[][cc][1.1]{\renewcommand{\arraystretch}{2}\begin{tabular}{c}{\tiny$t$}\end{tabular}}
   \captionsetup[subfloat]{captionskip=-1pt}
  \subfloat[\!Algorithm~\eqref{eq::Alg-EvnTrig}: $\alpha\!=\!1,~~~~$ $\!\beta\!=\!2$, $\!\Delta\!=\!0.2$~s.]{
    \includegraphics[height=.965in]{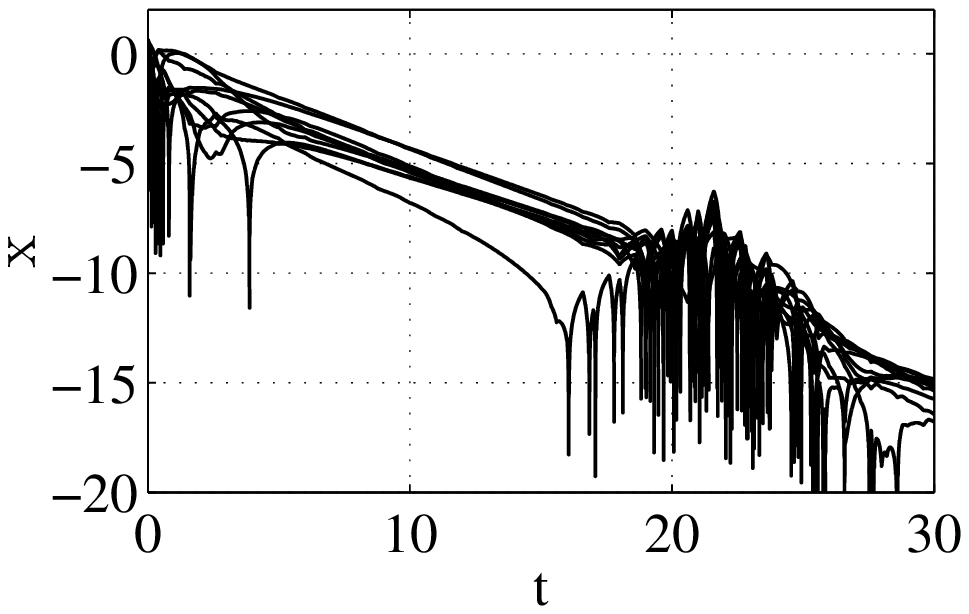} 
  }\quad
  \subfloat[\!Euler discretization of~\eqref{eq::Alg} ($\alpha\!=\!1$, $\beta\!=\!1$, $\!\Delta\!=\!0.2$~s.)]{
    \includegraphics[height=.965in]{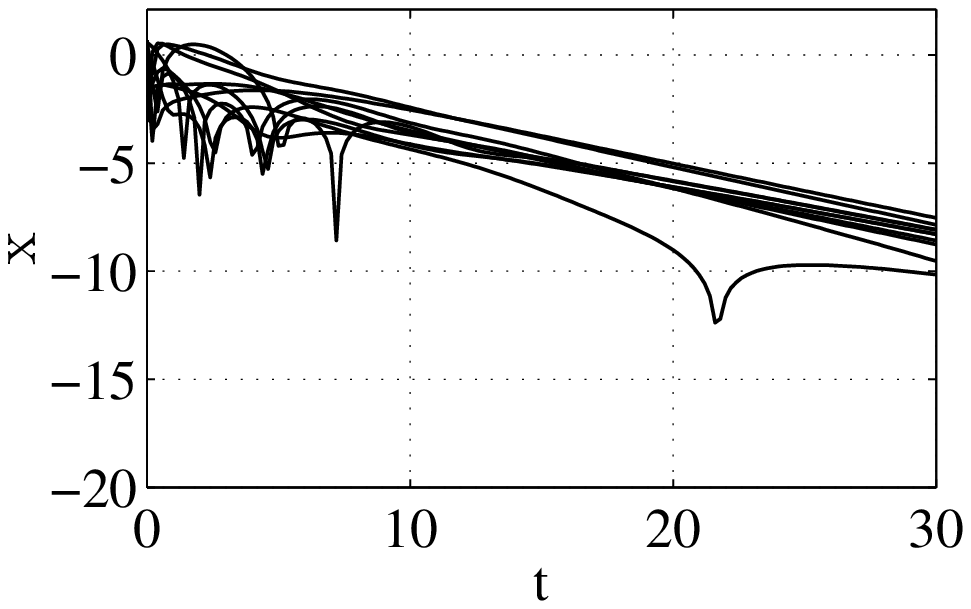} 
  }
  \vspace{-2pt}
    \caption{Performance evaluation of the
    algorithm~\eqref{eq::Alg-EvnTrig}  when the communication is
    periodic vs. the Euler-discretized implementation of the
    algorithm~\eqref{eq::Alg}.
  }\label{fig::Ex2sim-p}
\end{figure}

Figures~\ref{fig::Ex2sim-p}(a)-(b) compare the evolution of the agents
for~\eqref{eq::Alg-EvnTrig} with periodic communication and for an
Euler discretization of~\eqref{eq::Alg} over the network
in~Figure~\ref{fig::network}. In these simulations, we fixed $\Delta$
and varied $\beta$ until the algorithm becomes close to the
divergence. The results show that~\eqref{eq::Alg-EvnTrig} can use a
larger $\beta$, which reveals that, for the same amount of
communication effort,~\eqref{eq::Alg-EvnTrig} achieves faster
convergence.

\begin{figure}[t!]
  \unitlength=0.5in
  \psfrag*{x}[][cc][1.1]{\renewcommand{\arraystretch}{0.5}\begin{tabular}{c}$\ln(|x^i-\mathsf{x}^\star|)$\\~\end{tabular}}
  \psfrag*{A}[][cc][1.1]{\renewcommand{\arraystretch}{0.5}\begin{tabular}{c}Agents\\~\end{tabular}}
  \psfrag*{t}[][cc][1.1]{\renewcommand{\arraystretch}{2}\begin{tabular}{c}$t$\end{tabular}}
  \centering {
    \includegraphics[height=1.4in]{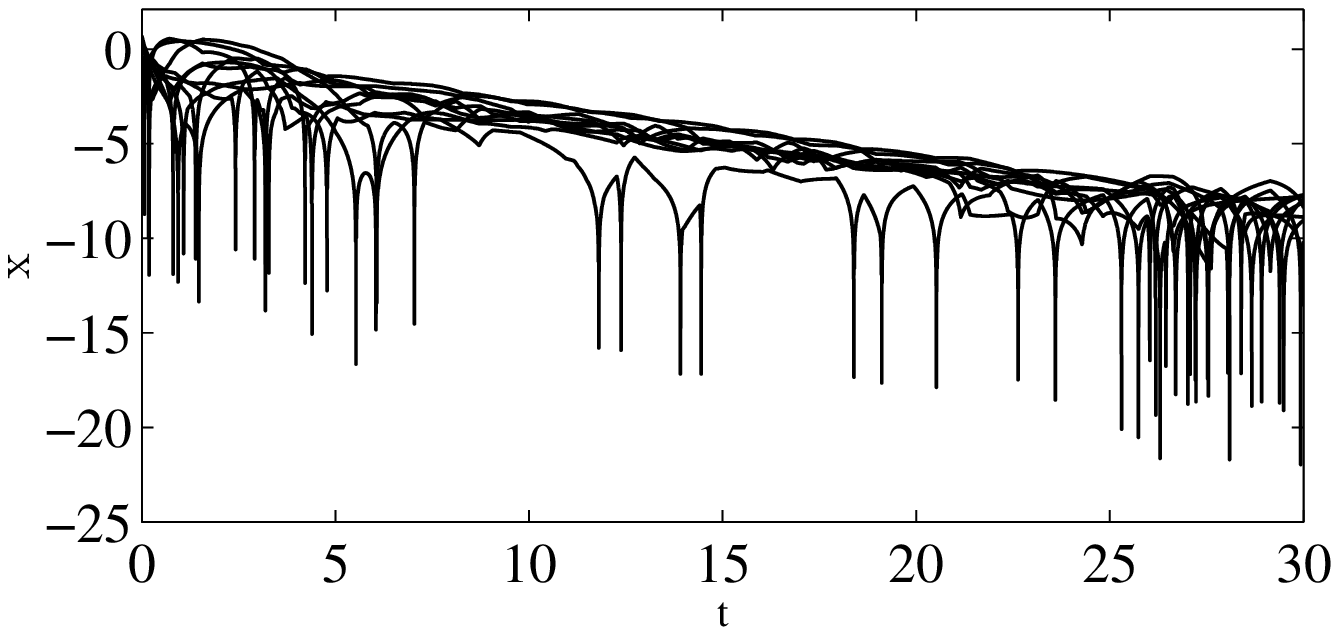}
  }
  \\\vspace{-5pt}
  \subfloat
  {
    \includegraphics[height=1.4in]{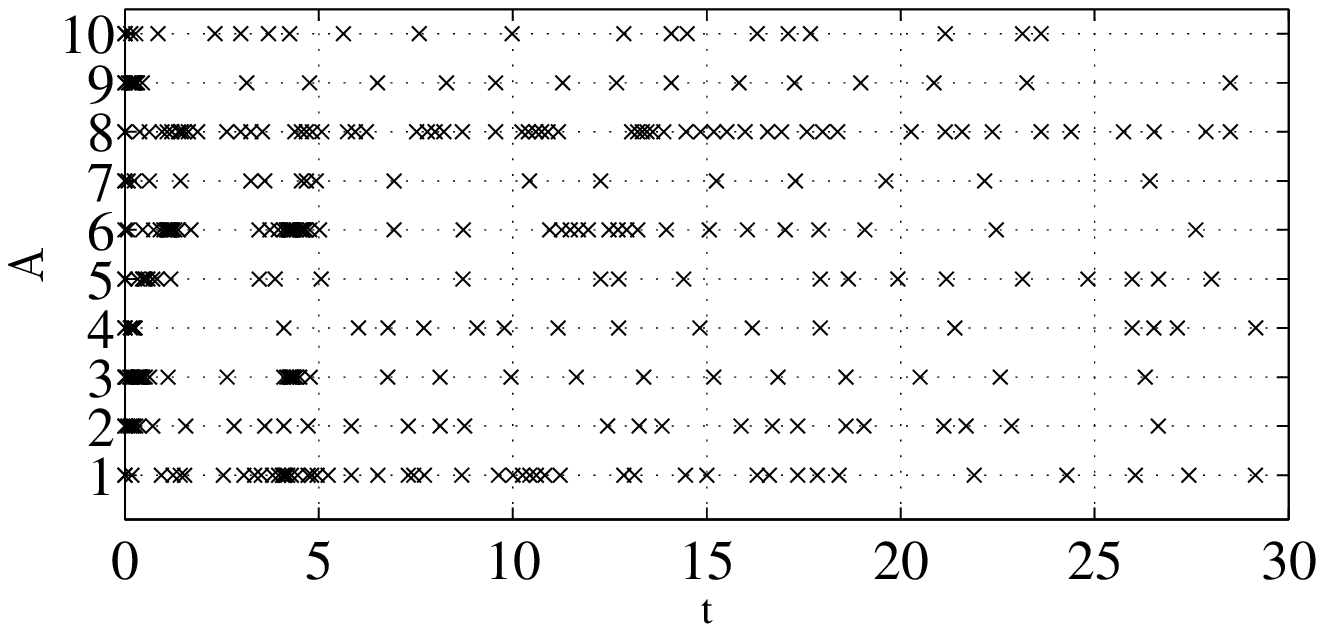}
  } \vspace{-2pt}
  \caption{Execution of
    algorithm~\eqref{eq::Alg-EvnTrig} using $\alpha=\beta=1$
    when the distributed event-triggered communication
    law~\eqref{eq::TrigLaw_Distributed} with $\eps^i=0.002$,
    $i\until{N}$ is employed: 
in the bottom plot $\times$ shows the time an event is triggered by an agent.}\label{fig::Ex3sim}
\end{figure}

Figure~\ref{fig::Ex3sim} shows the time history of the natural log
error $x^i-\mathsf{x}^\star$ and~the communication execution times of
each agent $i \until{10}$ of~\eqref{eq::Alg-EvnTrig} with the
distributed event-triggered communication
law~\eqref{eq::TrigLaw_Distributed}.
The results illustrate the behavior guaranteed by
Theorem~\ref{thm::UD_DecentralziedEventTriggered}: the communication
times are asynchronous, the operation is Zeno-free, and the states
converge to an $\|\vect{\eps}\|^2$-neighborhood of~$\sf{x}^\star$.

\section{Conclusions}\label{sec::conclu}

We have presented a novel class of distributed continuous-time
coordination algorithms that solve network optimization problems where
the objective function is strictly convex and equal to a sum of local
agent cost functions.  For strongly connected and weight-balanced
agent interactions, we have shown that our algorithms converge
exponentially to the solution of the optimization problem when the
local cost functions are strongly convex and their gradients are
globally Lipschitz. This property is preserved in dynamic networks as
long as the topology stays strongly connected and weight-balanced. For
connected and undirected agent interactions, we have shown that
exponential convergence still holds under the relaxed conditions of
strongly convex local cost functions with locally Lipschitz
gradients. In this case, asymptotic convergence also holds when the
local cost functions are just convex. We have also explored the
implementation of our algorithms with discrete-time
communication. Specifically, we have established asymptotic
convergence under periodic, centralized synchronous, and distributed
asynchronous event-triggered communication schemes, paying special
attention to establishing the Zeno-free nature of the algorithm
executions. Future work will focus on strengthening the results to eliminate the
offline computation of the design parameters, the study of the
robustness against disturbances, time delays, and asynchronous agents'
clocks, the exploration of agent abstractions for self-triggered
implementations, and the use of triggered control methods in other
coordination problems, including constrained, time-varying, and online
scenarios, and networked games.

{\footnotesize%

}

\end{document}